\documentclass[11pt]{article}

\usepackage{amsmath}
\usepackage{amssymb} 
\usepackage{caption}
\usepackage{graphicx}
\usepackage{authblk}
\usepackage[hypertexnames=false,colorlinks=true,linkcolor=blue,citecolor=blue]{hyperref}
\usepackage[numbers,comma,square,sort&compress]{natbib}
\usepackage[a4paper,text={6.5in,10in},centering]{geometry}

\graphicspath{{eps/}{pdf/}{images/}}
\makeatletter\@addtoreset{equation}{section}\makeatother

 {\begin{trivlist} \item[]{\bf Proof. }}%
 {\hspace*{\fill}$\rule{.4\baselineskip}{.4\baselineskip}$\end{trivlist}}

 {\begin{trivlist}\item[]\textbf{Acknowledgments.}}{\end{trivlist}}
{\begin{trivlist}\item[]\textbf{Proof#1 }}%
 {\hspace*{\fill}$\rule{0.3\baselineskip}{0.35\baselineskip}$\end{trivlist}}

\usepackage{amsmath}
\usepackage{caption}
\usepackage{subcaption}
\usepackage{wrapfig}
\def\Ab{\mathbf A}
\def\Xb{\mathbf X}
\def\xb{\mathbf x}
\def\Ub{\mathbf U}
\def\Vb{\mathbf V}

\def\vb{\mathbf v}
\def\Sigmab{\mathbf \Sigma}

\begin{document}

\title{Combining Dynamic Mode Decomposition with Ensemble Kalman Filtering for Tracking and Forecasting}
\author[1]{Stephen A Falconer}
\author[1]{David J.B. Lloyd}
\author[1]{Naratip Santitissadeekorn}
\affil[1]{\small Department of Mathematics, University of Surrey, Guildford, GU2 7XH, UK}
\date{\today}
\maketitle

%%%%%%%%%%%%%%%%%%%%%%%%%%%%%%%%%%%%%%%%%
\begin{abstract}

Data assimilation techniques, such as ensemble Kalman filtering, have been shown to be a highly effective and efficient way to combine noisy data with a mathematical model to track and forecast dynamical systems. However, when dealing with high-dimensional data, in many situations one does not have a model, so data assimilation techniques cannot be applied. In this paper, we use dynamic mode decomposition to generate a low-dimensional, linear model of a dynamical system directly from high-dimensional data, which is defined by temporal and spatial modes, that we can then use with data assimilation techniques such as the ensemble Kalman filter. We show how the dynamic mode decomposition can be combined with the ensemble Kalman filter (which we call the DMDEnKF) to iteratively update the current state and temporal modes as new data becomes available. We demonstrate that this approach is able to track time varying dynamical systems in synthetic examples, and experiment with the use of time-delay embeddings. We then apply the DMDEnKF to real world seasonal influenza-like illness data from the USA Centers for Disease Control and Prevention, and find that for short term forecasting, the DMDEnKF is comparable to the best mechanistic models in the ILINet competition.
\end{abstract}

\section*{Keywords}
Dynamic mode decomposition; Ensemble Kalman filter; Data-driven modelling; Data assimilation; Dynamical systems

\section{Introduction}

Data assimilation refers to the collection of methods that integrate vast data sets with sophisticated mathematical models, to track and forecast systems that may evolve or change \cite{data_assimilation_book}. The majority of its applications lie in the earth sciences \cite{da_use_earth_sciences}, however due to the generality of its techniques they have also been successfully applied in a wide range of areas from medicine \cite{da_use_medicine} to ecology \cite{da_use_ecology}. The Kalman filter \cite{kf_og} is one such data assimilation technique widely used throughout industry \cite{kf_applications} that optimally combines predictions from a linear model with Gaussian data. Whilst traditionally applied to a model's state, the parameters of the model can simultaneously be filtered, leading to what is known as the joint state-parameter estimation problem \cite{kf_state_param}. If the system being filtered is nonlinear, alternative versions of the Kalman filter can be utilized such as the extended Kalman filter \cite{ekf}, unscented Kalman filter \cite{ukf} or ensemble Kalman filter (EnKF) \cite{enkf}. The EnKF represents the distribution of a system's state with an ensemble of random samples, that can then be used to estimate useful statistics like the state's covariance via the sample covariance or a point estimate of the state via the sample mean \cite{enkf} and is well-suited for high-dimensional problems. All of these methods require a model of the system, however if no model exists then one must be generated and the most generalizable way to do this is via data-driven modelling.

Dynamic mode decomposition (DMD) is a data-driven modelling technique for identifying low dimensional, spatial and temporal patterns within a dynamical system directly from high-dimensional data \cite{official_intro}. It does this by postulating the state vector is evolved via a linear system and looking for a low-dimensional approximation of the eigenvalues (temporal modes) and corresponding eigenvectors (spatial modes). Spatial modes can be thought of as modes that decompose state variables into separate components that evolve together linearly in time. The corresponding temporal modes describe whether a spatial mode is growing, decaying or stationary in time. DMD has been used to approximate dynamical systems from measurement data in a multitude of fields, ranging from epidemiology \cite{epidemiology_use}, finance \cite{finance_use} and neuroscience \cite{neuroscience_use}. Due to its popularity, it has been extended to systems that are nonlinear in their recorded measurement functions via Extended/Kernel DMD \cite{extended_dmd}/\cite{kernel_dmd}, with one such extension Hankel-DMD \cite{hankel-dmd} employing time-delay embeddings of the original observables. In the presence of measurement noise, the standard DMD has been shown to induce a systematic bias by asymmetrically attributing all noise to the model's target output measurements and none to its inputs during training \cite{ncdmd}. This systematic bias, prompted the creation of noise handling variants of DMD that directly account for the noise term \cite{ncdmd}, the Forward Backward DMD \cite{ncdmd} that performs DMD forwards and backward in time and combines the results, and Total DMD (TDMD) \cite{tdmd} that minimizes the total least squares error as opposed to minimizing the ordinary least squares error.

The aim of this paper is to develop an algorithm that iteratively improves the temporal modes (eigenvalues) and state estimates produced by DMD with the EnKF as new data becomes available. This would be highly useful for dynamical systems that make a change from growing or decaying behaviour over time. While estimating just the state of the system using the DMD modes can be done using a standard Kalman filter, without also filtering the model's temporal mode, estimates are likely to suffer if the system changes over time. Methods already exist that combine DMD with the Kalman filter \cite{kfdmd} or extended Kalman filter \cite{ekfdmd}, which apply filtering to estimate the entire system dynamics matrix. The filtering in our work is instead focused on efficiently tracking the system’s temporal modes, and forecasting the system’s future states. DMD produces a linear model which makes it a natural fit for the Kalman filter, however when a system’s state and temporal modes are estimated simultaneously the filtering process becomes nonlinear. Hence, we need to use a filter designed for a nonlinear model, and we chose the EnKF due to its versatility, scalability to large dynamical systems, and ease of implementation \cite{enkf_flexible}. While any DMD variant that produces temporal modes would be compatible with the DMDEnKF framework, we use TDMD to remain consistent with the EnKF's assumption that noise is present in the data. In tandem, we apply the DMDEnKF using a total least squares version of Hankel-DMD, henceforth referred to as the Hankel-DMDEnKF, to investigate the effect time-delay embeddings have on our framework.

To demonstrate the DMDEnKF method, we first test it on synthetically generated datasets. Initially, on a simple noisy oscillating system with a decreasing period of oscillation, we use the DMDEnKF to track the system's temporal modes and compare results with the Hankel-DMDEnKF, other iterative DMD variants, and ``gold standard" filtering methods. Next, we simulate a pandemic and evaluate the DMDEnKF's ability to track the system's temporal modes and generate multistep ahead forecasts.

\begin{figure}[htbp]
\begin{subfigure}[b]{\textwidth}
         \includegraphics[width=\textwidth]{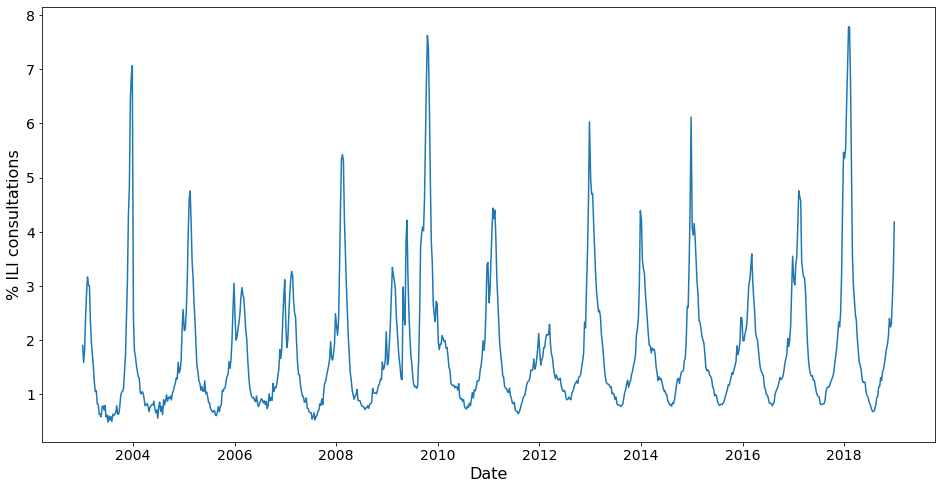}
     \end{subfigure}
     \caption{\centering\label{fig:ili_time_series}ILI consultations as a percentage of total weekly GP consultations in the US from 2003 to end of 2018. The data shows the annual peaks in ILI consultations that vary in size, timing and shape, which would make them difficult to model with a simple SIR-type model.}
\end{figure}

Finally, we apply the DMDEnKF and Hankel-DMDEnKF to real seasonal influenza-like illness (ILI) data in the United States from the Centers for Disease Control and Prevention (CDC) ILINet \cite{cdc_ili_definition} shown in Figure \ref{fig:ili_time_series}, with the aim of investigating their forecasting skills for ILI consultation rates. ILI is defined as a fever with a cough or sore throat that has no known cause other than influenza \cite{cdc_ili_definition} and infects between 9 and 35 million people in the US each year \cite{ili_deaths_updated}. Due to its prevalence, a multitude of methods have already been developed to model the spread of ILI \cite{ili_modelling_review1, ili_modelling_review2} and the approaches these models take can broadly be classified as either mechanistic or statistical \cite{ili_mech_or_stat}. Mechanistic methods \cite{ili_mechanistic1, ili_mechanistic2} make explicit hypotheses about what is driving the spread of an infectious disease, before then fitting parameters in the proposed models to the data. They have the advantage of being highly interpretable making them useful when trying to understand how one can control the spread of a disease \cite{ili_mechanistic_benefits}, however can make assumptions that are oversimplified \cite{mechanistic_vsml}. For example, a simple SIR-type model would struggle to describe specific behaviours like the drop in ILI consultations around Christmastime seen in Figure \ref{fig:ili_time_series}. Statistical methods \cite{ili_statistical1, ili_statistical2} are generally more versatile as they require fewer domain-specific assumptions, but both methods achieve a similar predictive skill in real time on the ILINet dataset \cite{ilinet_model_comparison}. The DMDEnKF attempts to find a middle ground between the two methods, remaining versatile by virtue of being purely data-driven but also providing some level of interpretability via the associated DMD modes.

The remainder of this paper will be structured as follows. First, a brief summary of DMD, Hankel-DMD and EnKF algorithms for completeness. After which, the DMDEnKF algorithm will be described in full. We will then apply the DMDEnKF and Hankel-DMDEnKF to synthetic data and compare their performance against other pre-existing, iterative DMD variants. Finally, we will use the DMDEnKF and Hankel-DMDEnKF on ILINet data to forecast the rate of ILI consultations up to 4 weeks into the future and examine their performance.

\section{DMDEnKF}

\subsection{Dynamic Mode Decomposition (DMD)}

Consider an $n$ dimensional state $\xb_k\in{\mathbb R^{n}}$ measured at regular time intervals $k=1,... ,m$. Assuming this time-series data was generated by a linear dynamical system, the consecutive states $\xb_k$ and $\xb_{k+1}$ are connected via \begin{equation}
    \xb_{k+1} = \Ab\xb_k 
    \label{eqn:axk=xk+1}
\end{equation}
for some unknown matrix $\Ab\in{\mathbb R^{n\times n}}$. By denoting \begin{equation}
    \Xb = \left[\begin{array}{cccc}
         |&|& &|  \\
         \xb_1&\xb_2&...&\xb_{m-1} \\
         |&|& &|
    \end{array}\right],\\
    \quad
    \Xb' = \left[\begin{array}{cccc}
         |&|& &|  \\
         \xb_2&\xb_3&...&\xb_{m} \\
         |&|& &|
    \end{array}\right],
    \label{eqn:XX'}
\end{equation}
equation \eqref{eqn:axk=xk+1} can be written succinctly over all consecutive data pairs as
\begin{equation}
    \Xb' = \Ab\Xb
    \label{eqn:AX=X'}.
\end{equation}
To minimize the mean squared error term $\sum_{k=1}^{m-1}\|\xb_{k+1}-\Ab\xb_{k} \|_2^2$, the standard DMD defines
\begin{equation}
    \Ab = \Xb'\Xb^+,
    \label{eqn:A=X'X+}
\end{equation}
where $\Xb^+$ is the Moore-Penrose pseudoinverse \cite{moore-pen} of $\Xb$. Efficiently solving for the eigendecomposition of $\Ab$ is the primary purpose of DMD, as these eigenvalues/eigenvectors correspond to spatio-temporal patterns in the data.

The DMD method starts by applying the Singular Value Decomposition (SVD) to the data matrix $\Xb$, representing it as the matrix multiplication of 2 real-valued, orthonormal matrices (complex and unitary if $\Xb \in \mathbb{C}^{n\times m}$) $\Ub \in \mathbb{R}^{n\times n},\Vb \in \mathbb{R}^{m\times m}$ and a rectangular diagonal matrix with decreasing non-negative real values ($\Sigmab \in \mathbb{R}^{n\times m}$) in the form
\begin{equation}
    \Xb = \Ub\Sigmab\Vb^*.
\label{eqn:X=USV}    
\end{equation}
The best rank $r$ approximation of a matrix according to the Eckart-Young Theorem \cite{eckart-young} is obtained by truncating its SVD, hence by truncating equation \eqref{eqn:X=USV} to a suitable rank $r$ \cite{optimal_trunc} we can compress the data matrix with minimal loss of information, which we write as
\begin{equation}
    \Xb \approx \Ub_r\Sigmab_r\Vb_r^*.
    \label{eqn:X=USVr}
\end{equation}
By performing this compression, we are implicitly assuming that there exists a low dimensional ($\leq r$), linear structure within the high-dimensional data.

The Moore-Penrose pseudoinverse can be found directly from the SVD computed in equation \eqref{eqn:X=USV} as $\Vb\Sigmab ^{-1}\Ub^*$. We use the rank $r$ truncated matrices from equation \eqref{eqn:X=USVr} for reasons of efficiency, setting \begin{equation}
\begin{split}
    \Ab &= \Xb'\Xb^+,\\
    &\approx\Xb'\Vb_r\Sigmab_r^{-1}\Ub_r^*.
    \end{split}
    \label{eqn:A=XVSU}
\end{equation}
This approximation of $\Ab$ now acts only on an $r$ dimensional subspace defined by Col$(\Ub_r)$. Hence, we can restrict $\Ab$ onto this $r$ dimensional subspace (representing the largest $r$ POD modes of $\Xb$) and denote the restricted $\Ab \in \mathbb{R}^{r\times r}$ as
\begin{equation}
\begin{split}
    \Tilde{\Ab} &= \Ub_r^{*}\Ab\Ub_r, \\
    & \approx \Ub_r^{*}\Xb'\Vb_r\Sigmab_r^{-1}.
    \end{split}
    \label{eqn:A=UXVS}
\end{equation}
We calculate the eigenvalues ($\lambda_i$), and corresponding eigenvectors ($\vb_i$)  of $\Tilde{\Ab}$, and define 
\begin{equation}
    \mathbf{\Lambda} = \left[\begin{array}{ccc}
         \lambda_1&0&0  \\
         0&\ddots&0 \\
         0&0&\lambda_r
    \end{array}\right], \quad
    \mathbf{W} = \left[\begin{array}{cccc}
         |&|& &|  \\
         \vb_1&\vb_2&...&\vb_r \\
         |&|& &|
    \end{array}\right].
    \label{eqn:LW}
\end{equation}

Reconstructing the eigenvalues/eigenvectors of the original operator $\Ab$ will provide insights into the structure of the system \cite{eigeninfo} and allow us to propagate it forward in time. The eigenvalues of $\Tilde{\Ab}$ ($\mathbf{\Lambda}$) can be shown to be equal to the eigenvalues of the original operator $\Ab$ \cite{tu_review}, however recovering the original eigenvectors is more involved and can be done using either projected or exact DMD. We use the exact DMD method introduced by Tu et al. \cite{tu_review} as it finds the exact DMD modes $(\mathbf{\Phi})$ for all eigenvectors with non-zero $\lambda_i$, where $\mathbf{\Phi}$ is defined as
    \begin{equation}
        \mathbf{\Phi} = \Xb'\Vb_{r}\Sigmab_r^{-1}\mathbf{W}.
        \label{P=XVSW}
\end{equation}
DMD modes with zero eigenvalues have no effect on the system's dynamics, so this restriction of exact DMD is of little consequence. This method finds $\Ab$ such that $\Ab\Xb = \Xb'$ exactly provided $r \geq$ rank($\Xb$) and $\Xb$ and $\Xb'$ are linearly consistent \cite{tu_review}. 

With $\mathbf{\Lambda}$ and $\mathbf{\Phi}$ in hand, we can construct a $r$ dimensional approximation of $\Ab$, however still need to find the initial phase and amplitude of each mode. The standard method \cite{official_intro} for computing this vector ($\mathbf{b}$) is to rewrite the initial state $\xb_1$ in a basis of the DMD modes via
\begin{equation}
    \mathbf{b} = \mathbf{\Phi}^{+}\xb_1.
    \label{eqn:b=Px}
\end{equation}
It is worth noting that there exist alternative methods for example \cite{opt_dmd,sparse_dmd} that focus on optimizing $\mathbf{b}$ over all data points with additional conditions.

To summarise, the final solution to the discrete system can be written as
\begin{equation}
    \xb_k = \mathbf{\Phi\Lambda}^{k}\mathbf{b}.
    \label{eqn:x=PLb}
\end{equation}

In the remainder of the paper, we call $\mathbf{\Lambda}$ the temporal modes and $\mathbf{\Phi}$ the spatial modes.

\subsubsection{Hankel-DMD}
Hankel-DMD first augments the original, measured state $\xb_k\in{\mathbb R^{n}}$, by appending to it measurements of the state at the previous $d-1$ time steps
\begin{equation}
    h(\xb_k) = \left[\begin{array}{c c c c c}
              {\xb_{k}}^T & {\xb_{k-1}}^T & \hdots &  {\xb_{k-(d-1)}}^T
             \end{array}\right]^T,
    \label{eqn:h(x)=xkxk-1xk-d-1}
\end{equation}
 to form a new state $h(\xb_k) \in{\mathbb R^{dn}}$. This is known as a time-delay embedding, and we refer to $d$ as the delay-embedding dimension. Taking time-delay embeddings, $h(\xb_k)$, to be our new states, matrices $\mathbf{X}$ and $\mathbf{X'}$ from equation \eqref{eqn:XX'} now become 
 \begin{equation}
    \Xb = \left[\begin{array}{ccc}
         \xb_d & \hdots & \xb_{m-1}  \\
         \vdots & \ddots & \vdots \\
         \xb_1 & \hdots & \xb_{m-d}
    \end{array}\right],\\
    \quad
    \Xb' = \left[\begin{array}{cccc}
         \xb_{d+1} & \hdots & \xb_{m}  \\
         \vdots & \ddots & \vdots \\
         \xb_2 & \hdots & \xb_{m-(d-1)}
    \end{array}\right].
    \label{eqn:hXhX'}
 \end{equation}
 With $X$ and $X'$ defined above, Hankel-DMD proceeds exactly as the standard DMD algorithm, generating eigenvalues $\mathbf{\Lambda}$, DMD modes $\mathbf{\Phi}$ and their initial states $\mathbf{b}$ as described above. The original system can be reconstructed/forecast for all time steps from $d$ onwards, by applying equation \eqref{eqn:x=PLb} and restricting the result to the first $n$ rows.

\subsubsection{Iterative DMD Variants}
There exists other variants of DMD that are designed to be applied iteratively, and in this paper we will compare these with the DMDEnKF in their ability to track a system's eigenvalues and make future state predictions. Streaming DMD \cite{streaming_dmd} is an adaption of the standard DMD algorithm to efficiently process new data as it becomes available, and the noise aware variant Streaming TDMD \cite{streaming_tdmd} is the first variant we wish to compare against. The  second  method  we  will use  for  comparison  is Windowed DMD \cite{online_dmd}, where the standard DMD described above is applied over a sliding window of the $w$ most recent data snapshots only. The final method we will be comparing against is Online DMD \cite{online_dmd}, specifically the variant of this algorithm that places an exponentially decaying weight $\rho$ on the importance of past measurements.

\subsection{Ensemble Kalman Filter (EnKF)}

Consider a discrete-time, nonlinear dynamical system with a stochastic perturbation
\begin{equation}
    \mathbf{x}_{k} = \mathbf{F}(\mathbf{x}_{k-1}) + \mathbf{w}_k, \quad
    \mathbf{w}_k \sim \mathcal{N}(\mathbf{0},\mathbf{Q}_k),
    \label{eqn:kf_prop_eq}
\end{equation}
where $\mathbf{F}$ is a nonlinear function $\mathbf{F}:\mathbb{R}^{n} \to \mathbb{R}^{n}$, $\xb_k \in \mathbb{R}^{n}$ is the system's state, $\mathbf{w}_k \in \mathbb{R}^{n}$ is a stochastic perturbation and $\mathcal{N}$ is the normal distribution with mean $\mathbf{0}$ and covariance matrix $\mathbf{Q}_k$. A measurement equation that relates what we observe to the true state of the system is given by
\begin{equation}
    \mathbf{y}_k = \mathbf{H}(\mathbf{x}_k) + \mathbf{v}_k, \quad
    \mathbf{v}_k \sim \mathcal{N}(\mathbf{0},\mathbf{R}_k),
    \label{eqn:kf_update_eq}
\end{equation}
where $\mathbf{H}:\mathbb{R}^{n} \to \mathbb{R}^{l}$ is the system's observation operator, $\mathbf{y_k} \in \mathbb{R}^{l}$ is an observation of the system, $\mathbf{v_k} \in \mathbb{R}^{l}$ is the noise in the observation and $\mathcal{N}$ is the normal distribution with mean $\mathbf{0}$ and covariance matrix $\mathbf{R}_k$. We focus on the instance relevant to our use case where $\mathbf{H}$ is linear, so can be represented by a matrix $\mathbf{H} \in \mathbb{R}^{l\times n}$.

In general, filtering methods aim to combine information from the state-transition model \eqref{eqn:kf_prop_eq} and observation model \eqref{eqn:kf_update_eq} to compute the conditional density $p(\mathbf{x_k} | \mathbf{Y_k})$, where $\mathbf{Y_k}=(\mathbf{y_1}, ..., \mathbf{y_k}$). The Kalman filter is the optimal filter if $\mathbf{F}$ and $\mathbf{H}$ are both linear and the stochastic perturbations are normal \cite{kf_og}. The EnKF was developed to deal with the filtering problem where either the linear or normal assumption (or both) is violated \cite{enkf}. It exploits the Kalman formulation to propagate an ensemble of the state into a region of high probability in such a way that the ensemble spread would be consistent with the linear and normal model.

To begin the EnKF algorithm, an initial ensemble of $N$ state estimates $\hat{\xb}^{(1)}_0\text{,..., } \hat{\xb}^{(N)}_0$ is required. If an ensemble is not available, one can be generated from initial state estimates  $\hat{\xb}_0$ and covariance matrix $\mathbf{P}_0$ by taking $N$ independent draws from $\mathcal{N}(\hat{\xb}_0,\mathbf{P}_0)$. 

\paragraph{Algorithm}

The EnKF then acts as follows \cite{enkf}:

    \textbf{Step 1:} Propagate forward in time each ensemble member using equation \eqref{eqn:kf_prop_eq} for $i = 1,...,N$ via
    \begin{equation}
         \hat{\xb}^{(i)}_{k|k-1} = \mathbf{F}(\hat{\xb}^{(i)}_{k-1|k-1}) + \mathbf{w}^{(i)}_k.
            \label{eqn:kf_k|k-1_eqs}
    \end{equation}
    The notation $\hat{\xb}^{(i)}_{k|k-1}$ denotes the state estimate at time $k$ of the $i$th ensemble member $\hat{\xb}^{(i)}_k$ using only information up to time $k-1$, and $\hat{\xb}^{(i)}_{k-1|k-1}$ represents the same ensemble member at time $k-1$ using information up to time $k-1$. Each $\mathbf{w}^{(i)}_k$ is independently drawn from $\mathcal{N}(\mathbf{0},\mathbf{Q}_k)$. The current covariance matrix can also now be estimated via the sample covariance of the ensemble, which we denote as $\mathbf{\hat{P}}_{k|k-1}$. This can then be used to estimate the Kalman Gain matrix $\mathbf{\hat{K}}_k$ as
        \begin{equation}
        \mathbf{\hat{K}}_k = \mathbf{\hat{P}}_{k|k-1}\mathbf{H}^T(\mathbf{H}\mathbf{\hat{P}}_{k|k-1}\mathbf{H}^T + \mathbf{R}_k)^{-1}.
        \label{eqn:kf_K}
    \end{equation}

    \textbf{Step 2:} Calculate the measurement innovation utilizing equation \eqref{eqn:kf_update_eq}.

    From measurement $\mathbf{y}_k$, we again use $i = 1,...,N$ and generate simulated measurements
       \begin{equation}
        \mathbf{y}^{(i)}_k =\mathbf{y}_k + \mathbf{v}^{(i)}_k
        \label{eqn:kf_sim_measures}
    \end{equation}
    where each $\mathbf{v}^{(i)}_k$ is an independent draw from $\mathcal{N}(\mathbf{0},\mathbf{R}_k)$. These simulated measurements $\mathbf{y}^{(i)}_k$ are combined with the ensemble members $\hat{\xb}^{(i)}_{k|k-1}$ from equation \eqref{eqn:kf_k|k-1_eqs} to define $N$ measurement innovations
           \begin{equation}
        \mathbf{e}^{(i)}_k  = \mathbf{y}^{(i)}_k - \mathbf{H} \hat{\xb}^{(i)}_{k|k-1}.
        \label{eqn:kf_innov}
    \end{equation}
    
    The $\mathbf{e}^{(i)}_k$ represent samples from the distribution of the distance of the model's prediction from the measured value.

    \textbf{Step 3:} Combine the model estimates in equation \eqref{eqn:kf_k|k-1_eqs} and measurement innovation of equation \eqref{eqn:kf_innov} via the estimated Kalman gain from \eqref{eqn:kf_K} to update each ensemble member's state estimate
\begin{equation}
\hat{\xb}^{(i)}_{k|k} = \hat{\xb}^{(i)}_{k|k-1} + \mathbf{\hat{K}}_k\mathbf{e}^{(i)}_k.
\end{equation}
We can generate a point estimate for the state $\hat{\xb}_k$ using the mean of the $N$ updated ensemble members. This process then repeats every time a new state measurement becomes available, with the updated ensemble from the previous data point becoming the initial ensemble for the new one.

We combine these 2 previously described techniques to form the DMDEnKF. This new, hybrid method uses DMD to generate a low dimensional model of a dynamical system that is then iteratively improved by the EnKF as new data emerges.
\subsection{DMDEnKF}

We now describe how we carry out filtering of the temporal modes and state of the system, while keeping the spatial modes found by one's chosen version of DMD on the ``spin-up" fixed. We note that once we allow the temporal modes to vary with the spatial modes being fixed, these are no longer eigenvalues/eigenvectors, and we then call them temporal modes. Consider an $n$ dimensional state $\xb_k\in{\mathbb R^{n}}$ measured at regular time intervals $k=1,... ,m$ and then measured iteratively at times $k=m+1,...$.
\paragraph{Algorithm}
\hspace{1ex}

\textbf{Step 1:}
Perform the chosen version of DMD on the dataset $\xb_1,...,\xb_m$, defining $\Xb, \Xb'$ as before in equation \eqref{eqn:XX'} to obtain the expression
\begin{equation}
    \xb_k = \mathbf{\Phi\Lambda}^{k}\mathbf{b},
    \label{eqn:x=PLb2}
\end{equation}
where 
\begin{equation}
    \mathbf{\Lambda} = \left[\begin{array}{ccc}
         \lambda_1&0&0  \\
         0&\ddots&0 \\
         0&0&\lambda_r
    \end{array}\right], \quad
    \mathbf{\Phi} = \left[\begin{array}{cccc}
         |&|& &|  \\
         \mathbf{d}_1&\mathbf{d}_2&...&\mathbf{d}_r \\
         |&|& &|
    \end{array}\right], \quad
        \mathbf{\mathbf{b}} = \left[\begin{array}{c}
         b_1  \\
         \vdots \\
         b_r
             \end{array}\right],
    \label{eqn:LPb}
\end{equation}
and defining $\lambda_i$, $\mathbf{d}_i$, $b_i$ as the $i$th temporal mode, DMD mode, initial condition triplet of the $r$ retained modes. This acts as a spin-up process to generate a model we can then filter using the EnKF.

\textbf{Step 2:} Define the matrices required for the EnKF's ensemble initialisation, propagation via equation \eqref{eqn:kf_prop_eq}, and measurement using equation \eqref{eqn:kf_update_eq}.

First, rewrite each of the $r$ temporal modes in polar coordinates as 
\begin{equation}
    \lambda_i = \tau_i e^{\theta_i\mathrm{i}},
    \label{eqn:l=ret}
\end{equation}
where $\tau_i \geq 0$, $0 \leq \theta_i < 2\pi$ and $\mathrm{i}^2 = -1$.
As $\xb_k\in{\mathbb R^{n}}$, the temporal modes in the DMD model's spectrum will either be real or in a complex conjugate pair. When filtering, we view the temporal modes as a time varying parameter. However, we must enforce that the real temporal modes remain real and complex conjugate pairs remain intact, as this ensures the state output by the model will still be real. We do this by defining the filterable model parameters $\mu_i$ as new variables for $i=1,... ,r$
\begin{equation}
  \mu_i=\left\{
  \begin{array}{@{}ll@{}}
    \tau_i, & \text{if}\ \theta_i = 0 \text{, or}\ \nexists\ j \text{ for } j < i \text{ such that } \lambda^*_j = \lambda_i,  \\
    \theta_i, & \text{otherwise}.
  \end{array}\right.
  \label{eqn:mu=tautheta}
\end{equation} 
Written in this way, these $\mu_i$'s represent all the possible degrees of freedom in the model's temporal modes under the additional constraint of producing a real state estimate. By maintaining a note of the positional indexes of each complex conjugate pair produced in the initial DMD, it is possible to recover the $\lambda_i$ representation from the $\mu_i$'s. While this transformation technically requires the full list of $\mu_i$'s, we informally write $\Lambda(\mu_i) = \lambda_i$ to symbolize the reversion from $\mu_i$'s back to $\lambda_i$'s.

\textbf{Ensemble initialisation:} We can now define an augmented joint parameter state $\mathbf{z}_0\in{\mathbb R^{n+r}}$ to be used as the initial state for the EnKF
\begin{equation}
            \mathbf{\mathbf{z}_0}
         = \left[\begin{array}{cccccc}
            - & \mathbf{x}_m & - &\mu_1 & \hdots & \mu_r
             \end{array}\right]^T.
    \label{eqn:z=xr}
\end{equation}
We denote the joint parameter state at time $m+k$ as $\mathbf{z}_k\in{\mathbb R^{n+r}}$. To generate an initial ensemble from this state, we first define sample covariance $\mathbf{C} = (1/m)(\Xb' - \mathbf{\Phi}\mathbf{\Lambda}\mathbf{\Phi}^{+}\Xb){(\Xb' - \mathbf{\Phi}\mathbf{\Lambda}\mathbf{\Phi}^{+}\Xb)}^T$, which represents the current state uncertainty based on prediction errors in the spin-up DMD. We then form the initial covariance matrix 
\begin{equation}
    \mathbf{P}_{0} = \left[\begin{array}{c|c}
         \mathbf{C}&\mathbf{0}  \\
         \hline
         \mathbf{0}& \alpha_2 I_{r}
    \end{array}\right],
\end{equation}
where $\alpha_2 > 0$, $I_r$ is the r-dimensional identity matrix, and the $\alpha_2I_r$ term determines the initial uncertainty in the spin-up DMD's temporal modes. Take independent draws from $\mathcal{N}(\mathbf{z}_0,\mathbf{P}_0)$ until the ensemble is sufficiently large. The optimal ensemble size will vary from problem to problem, adding ensemble members will increase accuracy but at the cost of computational efficiency.

\textbf{Propagation:} Using the notation $\mathbf{z}^i_k$ to signify the $i$th element of $\mathbf{z}_k$, we define the matrix $\mathbf{\Lambda}_{\mathbf{z}_k} \in{\mathbb R^{r\times r}}$ for state $\mathbf{z}_k$ as
\begin{equation}
    \mathbf{\Lambda}_{\mathbf{z}_k} = \left[\begin{array}{ccc}
         \Lambda(\mathbf{z}^{n+1}_{k})&0&0  \\
         0&\ddots&0 \\
         0&0&\Lambda(\mathbf{z}^{n+r}_{k})
    \end{array}\right].
    \label{eqn:lz=reiz}
\end{equation}
The EnKF's propagation equation can be written as
\begin{equation}
    \mathbf{z}_{k+1} = \left[\begin{array}{c|c}
         \mathbf{\Phi}\mathbf{\Lambda}_{\mathbf{z}_k}\mathbf{\Phi}^{+}&\mathbf{0}  \\
         \hline
         \mathbf{0}& I_{r}
    \end{array}\right]
    \mathbf{z}_k
    +
    \mathbf{w}_k.
    \label{eqn:zk+1=z+w}
\end{equation}
For convenience, we introduce notation $\mathbf{z}^{i:j}_{k} \in \mathbb{R}^{j-i+1}$ to denote the $i$th through to the $j$th element of $\mathbf{z}_k$ where $i \leq j$
\begin{equation}
    \mathbf{z}^{i:j}_{k} =
    \left[\begin{array}{c c c}
            \mathbf{z}^{i}_{k} & \hdots & \mathbf{z}^{j}_{k}
             \end{array}\right]^T
             \label{eqn:zki:j=zki...zkj}.
\end{equation}
Equation \eqref{eqn:zk+1=z+w} propagates $\mathbf{z}^{1:n}_{k}$ representing the state in the DMD framework $\xb_{m+k}$ forward in time using the standard DMD equation with the updated temporal modes from $\mathbf{\Lambda}_{\mathbf{z}_k}$. The vector $\mathbf{z}^{n+1:n+r}_{k}$ represents the current estimate of the temporal modes in their $\mu_i$ representation and is unchanged other than the addition of noise by the propagation equation, for although we assume the temporal modes vary in time no direction of drift in the parameters is explicitly foreknown. The vector $\mathbf{w}_k \in {\mathbb R^{n+r}}$ is a normally distributed variable $\mathbf{w}_k \sim \mathcal{N}(\mathbf{0},\mathbf{Q}_k)$, and this represents the uncertainty within the model of the system. We construct $\mathbf{Q}_k$ as follows,
\begin{equation}
    \mathbf{Q}_{k} = \left[\begin{array}{c|c}
         \alpha_1 I_{n}&\mathbf{0}  \\
         \hline
         \mathbf{0}& \alpha_2 I_{r}
    \end{array}\right],
    \label{eqn:R=I}
\end{equation}
where $\alpha_1$ and $\alpha_2$ are constants determined by the user such that $\alpha_2 \ll \alpha_1$. This construction with $\mathbf{Q}_k$ a diagonal matrix assumes model errors for each element of $\mathbf{z}_k$ are uncorrelated with one another. The condition $\alpha_2 \ll \alpha_1$ ensures that the state of the DMD system $\mathbf{z}^{1:n}_{k}$ changes significantly faster than its temporal modes $\mathbf{z}^{n+1:n+r}_{k}$, as parameters by definition should vary slowly in time. Furthermore, for the temporal mode's moduli being filtered, it prevents the strictly positive modulus dropping below 0.

\textbf{Measurement:} We write the EnKF's measurement equation as
\begin{equation}
    \mathbf{y}_{k} = \left[\begin{array}{c|c}
         I_{n}&\mathbf{0}  \\
         \hline
         \mathbf{0}& \mathbf{0}
    \end{array}\right]
    \mathbf{z}_k
    +
    \mathbf{v}_k,
    \label{eqn:yk=z+v}
\end{equation}
where $\mathbf{y}_k \in {\mathbb R^{n}}$ are observations of the DMD state $\xb_{m+k}$, and $\mathbf{v}_k \in {\mathbb R^{n}}$ is a normally distributed variable $\sim \mathcal{N}(\mathbf{0},\mathbf{R}_k)$ representing the noise in the measurements. We assume new measurements $\mathbf{y}_k$ to be available for the full DMD state $\mathbf{z}^{1:n}_{k}$ but not its temporal modes $\mathbf{z}^{n+1:n+r}_{k}$, as this is consistent with the format of the data used to generate the spin-up DMD model. We also assume uncorrelated measurement noise on each dimension of the state, so choose a diagonal matrix $\mathbf{R}_k$.

\textbf{Step 3}
State measurements $\xb_{m+k}$ at times $k=1,...$ are being iteratively generated. By setting $\mathbf{y}_k = \xb_{m+k}$ as each new measurement arrives, we can iteratively apply the EnKF to produce a hybrid estimate for $\mathbf{z}_k$ that combines model predictions from $\mathbf{z}_{k-1}$ and noisy measurement $\mathbf{y}_k$. A brief summary of how the EnKF does this is provided in Section 2.2, and a more expansive description can be found at \cite{enkf_tutorial}.

\textbf{Step 4:}
The state of the original system $\xb_{m+k}$ can be reconstructed from $\mathbf{z}_k$ by simply taking it's first $n$ elements $\mathbf{z}^{1:n}_{k}$. Predictions $p$ steps ahead at time $m+k$ can also be forecast from $\mathbf{z}_k$ via
\begin{equation}
    \xb_{m+k+p} = \mathbf{\Phi}\mathbf{\Lambda}_{\mathbf{z}_k}^p\mathbf{\Phi}^{+}\mathbf{z}^{1:n}_{k}.
    \label{eqn:x=PLzPz}
\end{equation}

The Hankel-DMDEnKF is defined algorithmically in exactly the same way, with the only difference being that Hankel-DMD is applied over the ``spin-up" period as opposed to standard DMD.

\section{Synthetic Applications}

\subsection{Comparison against other iterative DMD variants}
To test the DMDEnKF, we first apply it to data generated from a synthetic system with time varying eigenvalues, which we aim to track. The dynamics of this system are governed by the 2 dimensional rotation matrix, where the angle of rotation $\theta_k$ increases linearly from $\pi/64$ to $\pi/8$ over the course of 500 time steps. The evolution of the state $\xb_k$ of the system can hence be written as
\begin{equation}
    \xb_{k+1} = 
    \left[\begin{array}{cc}
         \cos({\theta_k})&-\sin({\theta_k})  \\
         \sin({\theta_k})&\cos({\theta_k})
    \end{array}\right]
    \xb_k,
    \quad
    \xb_1 =     
    \left[\begin{array}{c}
         1\\
         0
    \end{array}\right],
    \label{eqn:xk+1=rotx}
\end{equation}
where $\theta_k = \pi/64 + \frac{(k-1)(7\pi/64)}{499}$ and $k = (1,...,500)$.

\begin{figure}[!htbp]
\begin{subfigure}[b]{\textwidth}
         \centering
         \includegraphics[width=\textwidth]{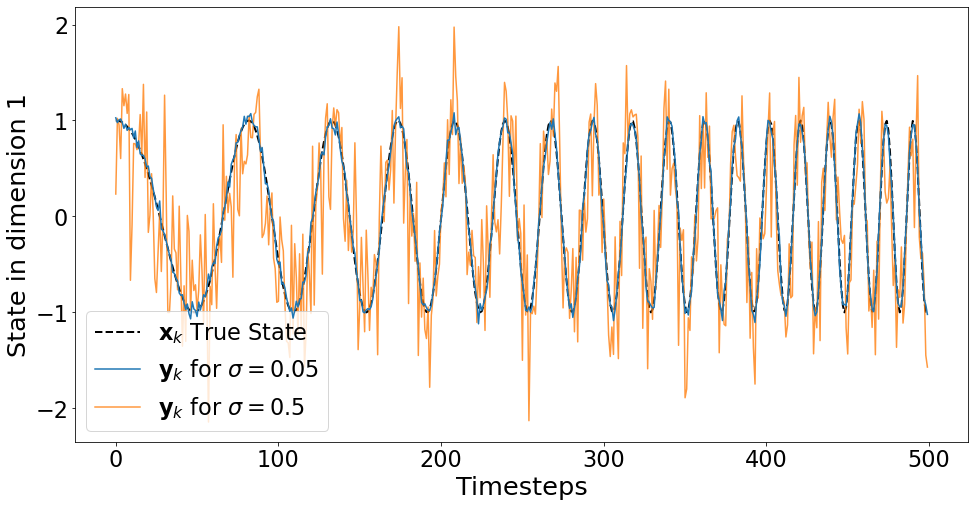}
     \end{subfigure}
     \caption{\centering\label{fig:synthetic_app_states}Time series for a synthetic system with a linearly increasing eigenvalue argument, showing the state's first dimension with no, low ($\sigma=0.05$) and high ($\sigma=0.5$) measurement noise.}
\end{figure}

We assume noisy measurement values $\mathbf{y}_k$ to be available for the state at each time step, such that
\begin{equation}
    \mathbf{y}_{k} = 
    \xb_k + \vb_k,
    \quad
    \vb_k \sim \mathcal{N}(\mathbf{0},\sigma^{2}I_{2}),
    \label{eqn:y=x+vvN0sI}
\end{equation}

 where each experiment $\sigma = 0.05 \text{ or } 0.5$ to simulate a low or high level of measurement noise respectively. The 500 values of $\mathbf{y}_k$ (shown in Figure \ref{fig:synthetic_app_states}) are used to train the DMDEnKF and Hankel-DMDEnKF, with the first 100 time steps being used for the spin-up process described in Step 1 of the DMDEnKF algorithm to produce the output described in equation \eqref{eqn:LPb}. 
 
 We will also train the iterative variants of DMD described at the end of Section 2.1 (Streaming TDMD\footnote{As the synthetic dataset is small, it is computationally tractable to apply batch methods to the data. Hence, instead of applying the true Streaming TDMD algorithm, we use batch TDMD over all data up to the current time step as a proxy for Streaming TDMD utilizing code from the PyDMD library \cite{pydmd}. As Streaming TDMD approximates the results of TDMD with the only differences occurring due to additional data compression steps in Streaming TDMD's algorithm, we believe this to be an acceptable substitution.}, Windowed DMD and Online DMD) on this dataset to compare their ability to track the system's time varying eigenvalues against that of the DMDEnKF. Within the Windowed DMD algorithm, we replace DMD with TDMD to allow for this method to effectively handle the noise in the data, henceforth referring to this amalgamation of the two methods as Windowed TDMD. To implement Online DMD, we use code made available by its creators here \cite{odmd_github}. Computational parameters were set as follows; window size $w = 10$ for Windowed TDMD, exponential decay rate $\rho = 0.9$ for Online DMD, delay-embedding dimension $d=50$ for the Hankel-DMDEnKF and spin-up time steps $m = 100$ for the DMDEnKF as previously stated.

At each time step $k$, the system's true eigenvalues can be written in modulus-argument form as
\begin{equation}
    \lambda_k = 1e^{\pm \theta_k i},
\end{equation}
and for each time step where the models are defined their estimates of the system's eigenvalues can also be written as
\begin{equation}
    \hat{\lambda}_k = \hat{\tau}_ke^{\pm \hat{\theta}_k i}.
\end{equation}
We start by comparing the errors in each method's estimate of the constant eigenvalue modulus ($\hat{\tau}_k - 1$). A thousand runs of the synthetic data were generated for each value of $\sigma$, and the difference of each method's eigenvalue modulus and argument from their true values at every time step after the spin-up period were collected. When any of the methods failed to identify the eigenvalues as a complex conjugate pair at a given time step in a run, the dominant eigenvalue's modulus was used for modulus error calculations. The average errors in the eigenvalue modulus estimates are shown in Table \ref{table:eig_mod_means}.

\begin{table}[]
\begin{center}
 \begin{tabular}{||c c c||} 
 \hline
 Iterative DMD & \multicolumn{2}{c||}{\centering Mean Eigenvalue Modulus Error} \\
 Variant & $\sigma =0.05$ & $\sigma=0.5$\\[0.5ex]
 \hline \hline
 Windowed TDMD & $9.82 \times 10^{-3}$ & 1.39 \\
 \hline
 Online DMD & $6.04 \times 10^{-3}$ & $3.06 \times 10^{-1}$ \\ 
 \hline
 Streaming TDMD  & $2.31 \times 10^{-4}$ & $2.50 \times 10^{-3}$ \\
 \hline \hline
 DMDEnKF & $8.07 \times 10^{-3}$ & $1.89 \times 10^{-2}$ \\
 \hline
 Hankel-DMDEnKF & $9.49 \times 10^{-3}$ & $1.38 \times 10^{-2}$ \\
 \hline
\end{tabular}
\end{center}
\caption{\centering Mean absolute errors in the synthetic system's eigenvalue modulus estimates produced by each iterative DMD variant over all time steps over the course of all 1000 runs. Measurement noise is set to either low levels with $\sigma = 0.05$ (left) or high levels with $\sigma = 0.5$ (right). Streaming TDMD scored significantly lower errors than all other methods, and as noise levels increased, errors in Windowed TDMD and Online DMD grew significantly larger than those produced by the DMDEnKF and Hankel-DMDEnKF.}
\label{table:eig_mod_means}
\end{table}
 
 For all levels of measurement noise, Streaming TDMD estimated the eigenvalue modulus the most accurately. This is due to the method's assumption of a stationary system, hence assigning an equal weight to the importance of each data point, which works well in the case of estimating a constant parameter. At low levels of measurement noise as seen in the first column of Table \ref{table:eig_mod_means}, Windowed TDMD, Online DMD, the DMDEnKF and Hankel-DMDEnKF all performed similarly well with mean eigenvalue modulus errors below 0.01. As errors in the eigenvalue modulus grow exponentially when forecasting future states, these 4 methods could produce acceptable short term forecasts but would quickly diverge from the true state as the forecast horizon was extended. At high levels of noise shown in the second column of Table \ref{table:eig_mod_means}, Windowed TDMD and Online DMD's eigenvalue modulus estimates degrade significantly, making them unsuitable for forecasting in this scenario. The errors in the DMDEnKF and Hankel-DMDEnKF remain fairly small, however are still an order of magnitude greater than those produced by Streaming TDMD.
 
\begin{figure}[h]
     
\begin{subfigure}[b]{0.48\textwidth}
         \centering
         \includegraphics[width=\textwidth]{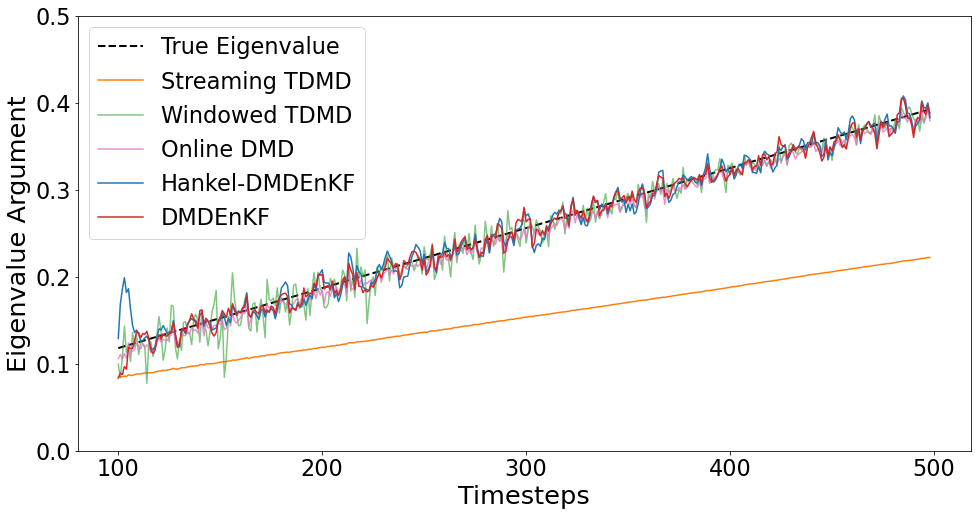}
                       \caption{\label{fig:arg_trajb}Eigenvalue argument trajectory for $\sigma = 0.05.$}
     \end{subfigure}
     \hfill
     \begin{subfigure}[b]{0.48\textwidth}
         \centering
         \includegraphics[width=\textwidth]{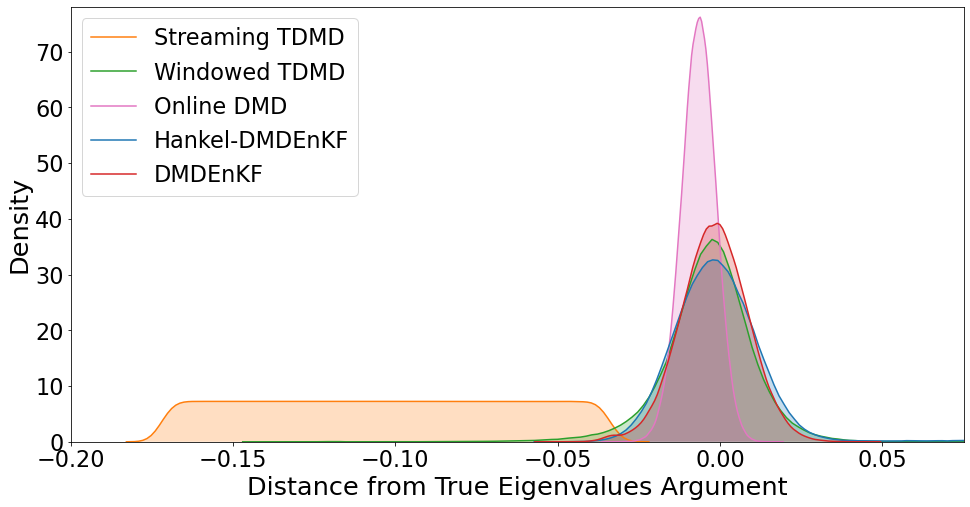}
         \caption{\label{fig:arg_distb}Error distribution for $\sigma = 0.05.$}
     \end{subfigure}
     \hfill
\begin{subfigure}[b]{0.48\textwidth}
         \centering
         \includegraphics[width=\textwidth]{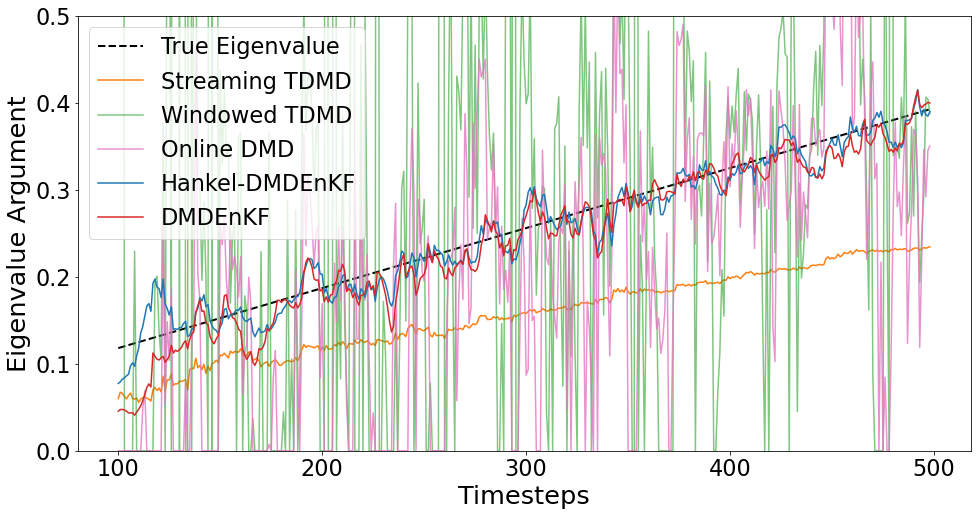}
                       \caption{\label{fig:arg_trajc}Eigenvalue argument trajectory for $\sigma = 0.5.$}
     \end{subfigure}
     \hfill
     \begin{subfigure}[b]{0.48\textwidth}
         \centering
         \includegraphics[width=\textwidth]{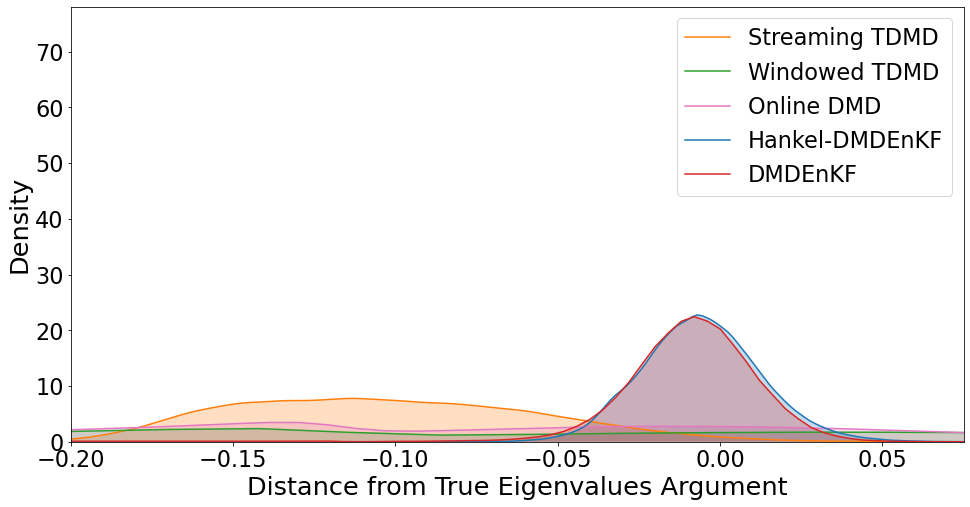}
         \caption{\label{fig:arg_distc}Error distribution for $\sigma = 0.5.$}
     \end{subfigure}
          \caption{\centering\label{fig:arg_traj} Estimates of the synthetic system's eigenvalue argument produced by each iterative DMD variant. Presented are typical trajectories of the eigenvalue argument at each time step over the course of 1 experiment's run (left) and error distributions of the difference between the true system's eigenvalue argument and the estimated eigenvalue argument over all time steps over the course of all 10 runs (right). Measurement noise is set to either low levels with $\sigma = 0.05$ (top) or high levels with $\sigma = 0.5$ (bottom). The DMDEnKF and Hankel-DMDEnKF experience similar errors to Online DMD and Windowed TDMD at low measurement noise, but track the eigenvalue argument much more accurately than them for high measurement noise.}
\end{figure}

 A typical trajectory of the eigenvalue argument estimates ($\hat{\theta}_k$) for each method over the course of one run from the end of the spin-up period onwards can be seen in Figures \ref{fig:arg_trajb} and \ref{fig:arg_trajc}. The error distributions for each method's eigenvalue argument estimates ($\hat{\theta}_k - \theta_k$) over all 1000 runs are plotted in Figures \ref{fig:arg_distb} and \ref{fig:arg_distc}.
 
At low levels of noise as seen in Figures \ref{fig:arg_trajb} and \ref{fig:arg_distb}, all 5 methods on average underestimated the eigenvalue argument of the system. This is to be expected as the eigenvalue argument is increasing with time, meaning that all but the last data pair available to each method would have been generated using an argument smaller than its current value. Streaming TDMD exhibited the worst performance, again due to its equal weighting of every data point, however in this instance being a negative quality as it hampers the model's ability to adapt to fresh data that reflects the changing parameter. Windowed TDMD, Online DMD, the DMDEnKF and Hankel-DMDEnKF all performed similarly. Online DMD produced a tighter error distribution, but with a slightly larger bias than Windowed TDMD. This suggests that Online DMD's soft thresholding reduces the model volatility caused by measurement noise compared to the hard cut-off employed by Windowed TDMD. For this same reason however, Online DMD is slower to adapt to new measurements than Windowed TDMD, leading to a larger bias below the system's true eigenvalue argument. The DMDEnKF and Hankel-DMDEnKF performed very similar to Windowed TDMD at this noise level, however tweaks to the magnitude of the DMDEnKF's system uncertainty matrix can be made to balance the speed of model innovation with its volatility and produce distributions closer to that of Online DMD if required.

At higher noise levels shown in Figures \ref{fig:arg_trajc} and \ref{fig:arg_distc}, the performance of Windowed TDMD and Online DMD significantly degrades. Placing a larger weight on more recent samples allowed these methods to quickly adapt to changes in the system's parameters, however as the noise increases this induces an extreme volatility in their respective models. The performance of Streaming TDMD is not largely changed from the low noise case, still lagging behind the true system values but somewhat insulated from the noise by its symmetric treatment of all data points. Here the benefit of explicit inclusion of measurement noise in the DMDEnKF framework becomes apparent, as at this noise level the DMDEnKF and Hankel-DMDEnKF are the only techniques tested capable of producing an accurate eigenvalue argument estimate.

Furthermore, here we see the first significant difference in the performance of the DMDEnKF and Hankel-DMDEnKF, as the DMDEnKF's error distribution has a thin tail extending down to $-\pi/8$ which is not present in the error distribution of the Hankel-DMDEnKF. These additional errors are caused by the spin-up of the DMD for the DMDEnKF method occasionally failing to identify the system's eigenvalues as a complex conjugate pair (empirically, this happens $\sim 3\%$ of the time), due to the increased noise in the data. When this happens, the DMDEnKF catastrophically fails for the EnKF is unable to generate complex eigenvalues from real ones regardless of how many future time steps are filtered due to its formulation in equation \eqref{eqn:mu=tautheta}. This failure of the DMDEnKF can be mitigated in the following way. If the errors produced by the DMDEnKF during the filtering stage are deemed too large (e.g. exceed a given threshold) for a prolonged period of time, then the spin-up DMD process can be rerun on an extended dataset consisting of the original spin up data, plus the newly available data used so far in the filtering step. By including more data in the spin-up process, the spin-up DMD model is more likely to successfully capture the signal component in the data as a pose to measurement noise, and hence produce eigenvalues with the same structure as those of the true system. Time-delay embeddings make the SVD step in the DMD algorithm more robust to measurement noise \cite{time_delay_pca}. Hence, while the Hankel-DMDEnKF is similarly restricted by the eigenvalues it can produce at the filtering stage, in all 1000 runs of the synthetic data the spin up Hankel-DMD was able to identify the system's eigenvalues to be a complex conjugate pair, so this was not an issue for the Hankel-DMDEnKF.

\subsection{Comparing against DMD with a particle filter}

Having compared the performance of the DMDEnKF against other iterative DMD variants, we now focus on evaluating the filtering component of the algorithm. Since the linear DMD model acts nonlinearly in the filter when applied to both the model's state and eigenvalues, we compare the EnKF filter with a particle filter. Particle filters \cite{particle_filter_og} have been show to converge to the optimal filter as the number of particles tends to infinity for general nonlinear models with non-Gaussian noise \cite{particle_filter_convergence}. However, particle filters are restricted to low dimensional systems only, as the number of particles required scales approximately exponentially with the dimension of the state \cite{particle_filter_high_dims}. Hence, we compare the DMDEnKF and Hankel-DMDEnKF with a DMD plus particle filter which we will take to be the ``gold standard" estimation to assess how well the EnKF does with the nonlinear filtering problem. 

We use the same synthetic system \eqref{eqn:xk+1=rotx} with a linearly increasing eigenvalue argument as in the previous subsection to generate data with high levels of measurement noise ($\sigma = 0.5$); a trajectory of which can be seen in Figure \ref{fig:synthetic_app_states}. Again, the time-delay embedding dimension $d=50$ for the Hankel-DMDEnKF, and the first 100 time steps are used to train a spin-up DMD model, with the next 400 used to filter the state and spin-up model's eigenvalues. 

The DMDEnKF's filter state thus has dimension 4 (2 state dimensions and 2 temporal modes), while the Hankel-DMDEnKF's filter state is of dimension 102 (100 state dimensions and 2 temporal modes). To generate a ``gold standard" solution, at the filtering step we use a particle filter with 10,000 particles, applying multinomial importance resampling \cite{particle_filter_resampling_schemes} every time the effective sample size falls below half the number of particles to avoid sample degeneracy \cite{particle_filter_0.5_effective_sample_size}. For the DMDEnKF and Hankel-DMDEnKF at the filtering step, we run the EnKF with varying numbers of ensemble members ($N$), to see if as $N$ increases their estimates mean and covariance will tend to that of the particle filter ensemble. We generated 1000 runs of the synthetic data to apply the particle filter/EnKF with each value of $N$ to and collected the errors in the eigenvalue argument estimates for each method at every time step.

 \begin{figure}[htbp]
     
     \begin{subfigure}[b]{0.69\textwidth}
         \centering
         \includegraphics[width=\textwidth]{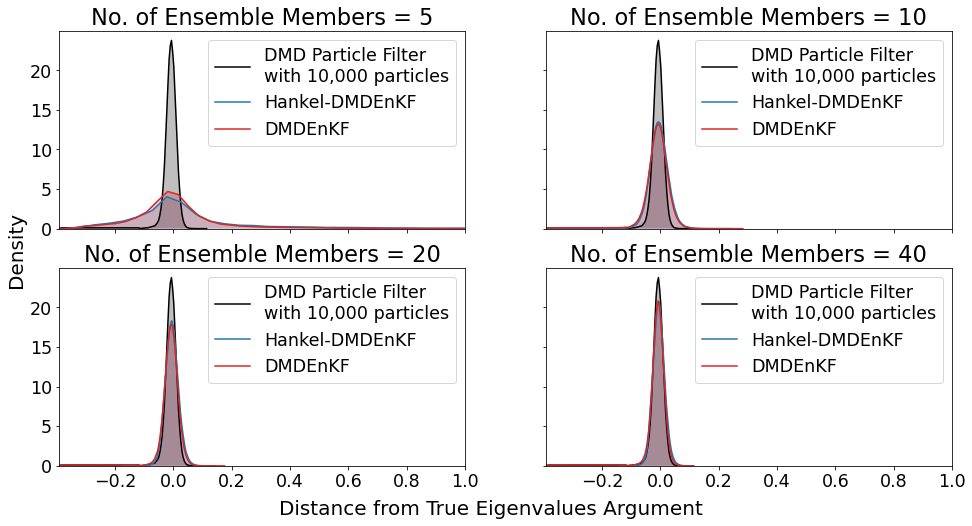}
         \caption{\label{fig:enkf_particle}Error distributions.}
     \end{subfigure}
          \hfill
               \begin{subfigure}[b]{0.3\textwidth}
         \centering
         \includegraphics[width=\textwidth]{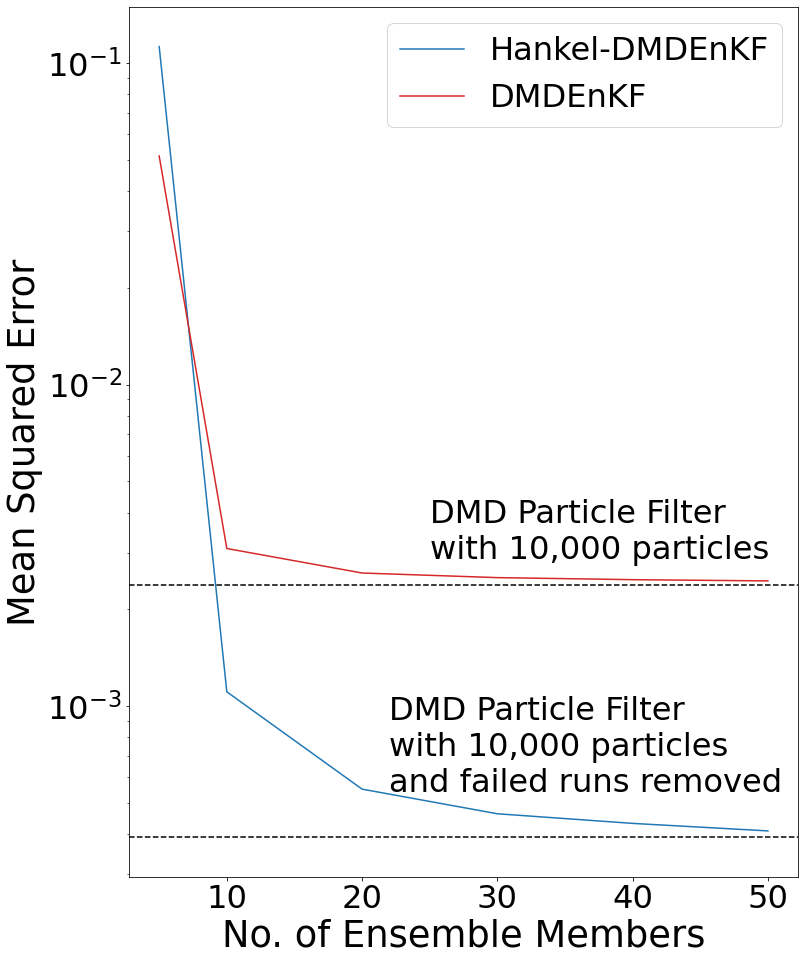}
         \caption{\label{fig:enkf_particle_mse}Mean squared errors.}
     \end{subfigure}
     \caption{\centering Error distributions (left) and mean squared errors (right) for estimates of the synthetic system's eigenvalue arguments produced by the DMDEnKF and Hankel-DMDEnKF with varying numbers of ensemble members ($N$) against those produced by a particle filter with 10,000 particles. Increasing $N$ quickly leads to error levels in the DMDEnKF and Hankel-DMDEnKF that are similar to those produced by their respective ``gold standards".}
\end{figure}

As can be seen in Figure \ref{fig:enkf_particle}, the DMD particle filter with 10,000 particles produces an extremely tight error distribution that is slightly biased to produce estimates below that of the true eigenvalue's argument. This is to be expected, as mentioned in the previous subsection, due to the system's eigenvalue argument constantly increasing. There is also a thin tail in the error distribution that extends down to $-\pi/8$. This is again a result of the spin up DMD sometimes failing to identify a complex conjugate eigenvalue pair, trapping the particle filter in the faulty assumption that the eigenvalues are real.

For low numbers of ensemble members ($N = 5$), the DMDEnKF and Hankel-DMDEnKF are centred at a similar value to the ``gold standard". However, they produce a far larger spread with long tails in both directions that imply a lack of robustness with this few ensemble members. With only a small increase to $N = 10$, both methods become more stable, as although they still have a larger variance than the particle filter, the long positive tails from $N=5$ have been eliminated. A similar pattern occurs as we move to $N=20$, with more ensemble members resulting in a tighter error distribution. At this point, the Hankel-DMDEnKF's distribution can be distinguished from that of the DMDEnKF and DMD particle filter by its aforementioned lack of a persistent thin negative tail. By $N=40$, the main peaks of the DMDEnKF, Hankel-DMDEnKF and ``gold standard" are almost indistinguishable on the graphical scale, with the DMDEnKF and DMD particle filter both sharing a thin negative tail.

Figure \ref{fig:enkf_particle_mse} shows how the mean squared error for the eigenvalue arguments predicted by the DMDEnKF and Hankel-DMDEnKF are affected by varying the number of ensemble members. For the DMDEnKF, errors initially sharply  decline as $N$ is increased, however on this small synthetic system returns diminish quickly after $N = 20$. By $N=50$, we achieve a mean squared error with the DMDEnKF only $\sim 3\%$ larger than that of the ``gold standard", despite using $200$ times fewer particles. When comparing the Hankel-DMDEnKF to the ``gold standard", the errors in the DMD particle filter's eigenvalue estimates are skewed by the runs in which the spin up DMD was unable to identify a complex conjugate eigenvalue pair, as Hankel-DMD did not encounter this problem on these synthetic examples. To attempt to fairly compare the filtering methods, we remove all runs in which the spin up DMD failed in this way, before again calculating the mean squared error for the DMD particle filter and recording it in Figure \ref{fig:enkf_particle_mse}. A similar pattern of reducing errors with diminishing returns can be seen for the Hankel-DMDEnKF as ensemble size is increased, and by $N=50$ its mean squared error is within $5\%$ of the newly calculated DMD particle filter's score. Our results show that in this simple, synthetic case at least, the EnKF is an efficient and effective solution to the nonlinear filtering problem that arise within the DMDEnKF framework.

\subsection{Tracking a synthetically generated pandemic}

\begin{figure}[htbp]
     
     \begin{subfigure}[b]{\textwidth}
         \centering
         \includegraphics[width=\textwidth]{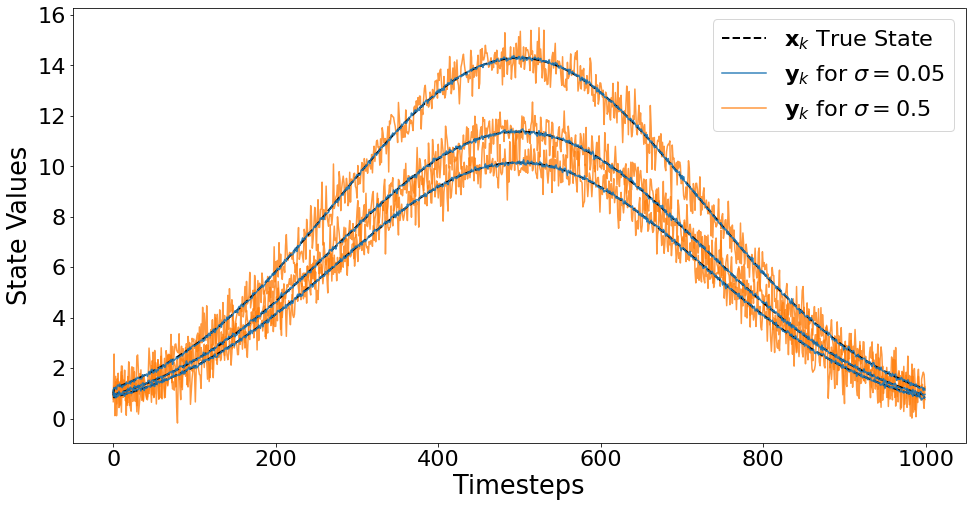}  \end{subfigure}
          \caption{\centering\label{fig:growth_decay_data} Time series for a synthetic system that simulates a pandemic, showing all 3 dimensions of the state with no, low ($\sigma=0.05$) and high ($\sigma=0.5$) measurement noise.}
\end{figure}

Lastly, we test the DMDEnKF's performance on synthetic data designed to simulate a simple pandemic with a state $\mathbf{x}_k$ representing the level of infection in 3 different population classes. The system's dynamics are governed by a matrix $\mathbf{A} \in \mathbb{R}^{3\times 3}$ that we randomly generate with non-negative elements each being drawn from the Uniform distribution $\mathcal{U}[0,1)$. The $(i,j)\text{th}$ element of $\mathbf{A}$ represents how the level of infection in class $j$ at time $k$ will affect the level of infection in class $i$ at time $k+1$. To control whether the synthetic pandemic is spreading or dying off, we then define a new matrix $\mathbf{\hat{A}} = \frac{\mathbf{A}}{\lambda_1}$ where $\lambda_1$ is the largest eigenvalue of $\mathbf{A}$, thus ensuring the spectral radius $\rho(\mathbf{\hat{A}}) = 1$. By introducing a constant $\gamma$, we can replace $\mathbf{A}$ with $\gamma \mathbf{\hat{A}}$ causing the state to grow if $\gamma > 1$ or decay for $\gamma < 1$. To simulate a pandemic, we linearly decrease $\gamma$ from 1.01 to 0.99 over the course of the experiment's 1000 time steps. The initial state used is a vector of ones. The system that generates the synthetic data can be written as
\begin{equation}
    \xb_{k+1} = 
    \gamma_k\mathbf{\hat{A}}
    \xb_k,
    \quad
    \xb_1 =     
    \left[\begin{array}{ccc}
         1 & 1 & 1
    \end{array}\right]^T,
\end{equation}
where the state $\xb_k \in \mathbb{R}^3$, $\gamma_k = 1.01 - \frac{0.02(k-1)}{999}$ and $k = (1,...,1000)$. We assume not to have access to the true state of the system $\mathbf{x}_k$ but instead noisy measurements $\mathbf{y}_k$ defined by
\begin{equation}
    \mathbf{y}_{k} = 
    \xb_k + \vb_k,
    \quad
    % \vb_k \sim \mathcal{N}(\mathbf{0},\sigma^{klukioijmhukjm 2}I_{3}).
\end{equation}
The constant $\sigma$ that governs the level of measurement noise is set to $\sigma = 0.05$ to represent low noise and $\sigma = 0.5$ for high noise as in \eqref{eqn:y=x+vvN0sI}. Figure \ref{fig:growth_decay_data} shows the values of the system's three state dimensions and the respective available measurements over the course of one run.

All five DMD variants tested had their computational parameters set to the same values as those used in the synthetic experiments in Section 3.1. The only small difference was that Streaming TDMD, Windowed TDMD, the DMDEnKF and Hankel-DMDEnKF truncated the data by removing the smallest singular value to reduce model instability caused by what was often a very low signal-to-noise ratio in this direction. Online DMD did not apply any truncation to the data as the method was not designed to do so, however it did not appear to suffer from any stability issues as a consequence.

\begin{figure}[htbp]
     
         \centering
         \includegraphics[width=\textwidth]{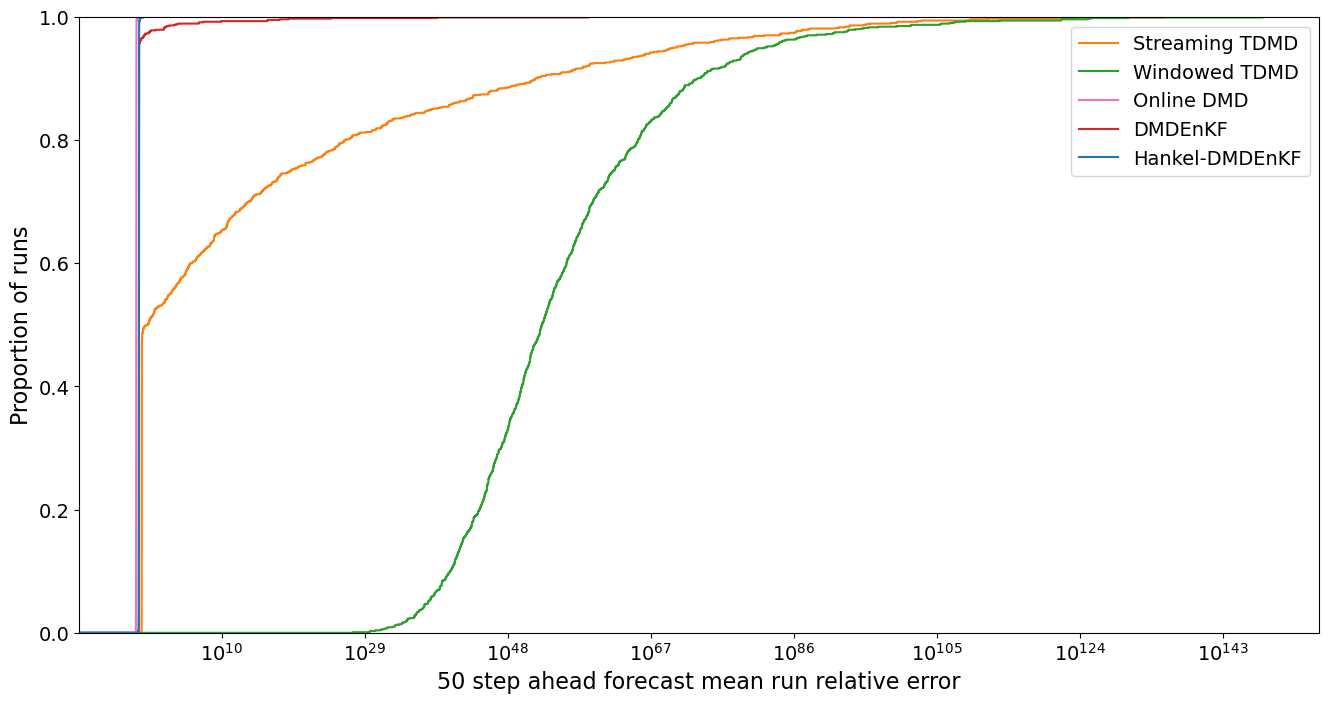}
                \caption{\centering\label{fig:0.05_error_distributions_all}
                Cumulative error distributions of the mean run relative errors for the 50 step ahead forecasts of each iterative DMD variant. Mean relative errors were calculated over all time steps for each run of the experiment, with the results from 1000 runs under low levels of measurement noise ($\sigma = 0.05$) displayed. Forecast errors had a wide range for some methods, due to exponentially compounding errors caused by forecasting 50 steps ahead. The DMDEnKF, Hankel-DMDEnKF and Online DMD produced errors orders of magnitude smaller than those of Streaming TDMD and Windowed TDMD.}
\end{figure}

The first 100 measurements $(\mathbf{y}_1,...,\mathbf{y}_{100})$ were used to initialize the models, and as each new data point $(\mathbf{y}_{100},...,\mathbf{y}_{1000})$ was successively fed into the models, they produced 50 step ahead forecasts $(\mathbf{\hat{x}}_{150},...,\mathbf{\hat{x}}_{1050})$. We generate 1000 data points, however a standard flu season lasts around 20 weeks. For this reason, we chose to forecast 50 steps ahead to mimic forecasting 1 week ahead in a more realistic timescale. The relative prediction errors $\hat{e}_{k} = \frac{\|\mathbf{x}_{k} - \mathbf{\hat{x}}_{k}\|}{\|\mathbf{x}_{k}\|}$ could then be calculated for $k=(150,...,1000)$ and the mean of these errors was the main metric we used to evaluate the forecasting skill of each method over the course of one run. A thousand runs were performed for both low and high levels of noise and the empirical cumulative distributions of 50 step ahead forecast mean run relative errors for low noise $(\sigma = 0.05)$ can be seen in Figure \ref{fig:0.05_error_distributions_all}.

The first noteworthy feature of the cumulative error distributions is the wide range in some method's forecast errors. This is a result of the 50 step ahead forecasts being produced by training each model to forecast 1 step ahead, then applying the trained model to the data 50 times. As such, forecast errors compound exponentially and small errors over a 1-step forecast horizon can become vast after 50 iterations. Inspecting the individual methods, we see Windowed TDMD to be the worst performing method. This is due to its aforementioned instability under measurement noise caused by considering only a small subset of the data at a time. This instability could be reduced by increasing the window size $(w)$ computational parameter, however as $w$ increases the model's ability to track a system that changes with time diminishes. Streaming TDMD had the second-largest errors, caused by the method's assumption of a stationary system hindering its ability to correctly adapt to the system's changing eigenvalues as new data became available. In the majority of cases, Online DMD, the DMDEnKF and Hankel-DMDEnKF all performed similarly well. All three methods exhibited cumulative error distributions tightly concentrated around a low error value, however in a few runs, the DMDEnKF became unstable and produced large errors. It is clear even at low levels of noise that the forecasting performance of Online DMD, the DMDEnKF and Hankel-DMDEnKF are far superior on this type of system to those of Windowed TDMD and Streaming TDMD. Hence, we now focus exclusively on these top three performing methods to allow for a thorough comparison of them on an appropriate error scale.

\begin{figure}[htbp]
\begin{subfigure}[b]{.48\textwidth}
         \includegraphics[width=\textwidth]{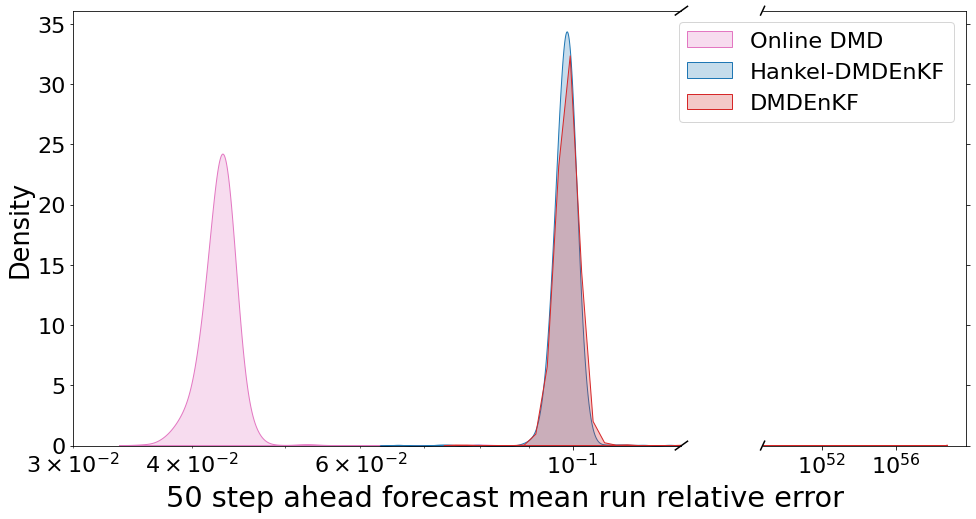}
                \caption{\centering\label{fig:0.05_error_distributions_odmd_dmdenkf}$\sigma = 0.05$}
\end{subfigure}
\hfill     
\begin{subfigure}[b]{.48\textwidth}
         \includegraphics[width=\textwidth]{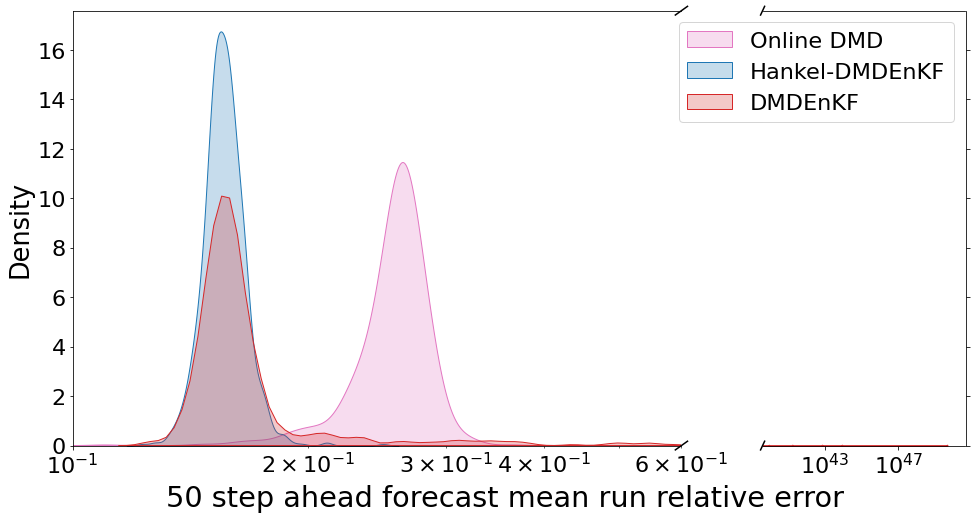}
                \caption{\centering\label{fig:0.5_error_distributions_odmd_dmdenkf}$\sigma = 0.5$}
\end{subfigure}
     \caption{\centering\label{fig:error_distributions_odmd_dmdenkf}Error distributions of the mean run relative errors for the 50 step ahead forecasts of Online DMD, the DMDEnKF and Hankel-DMDEnKF, attained over 1000 runs under low $\sigma=0.05$ (left) and high $\sigma=0.5$ (right) levels of measurement noise. Similar errors were found at both noise levels, with Online DMD performing better at low measurement noise and the DMDEnKF/Hankel-DMDEnKF performing better at high measurement noise.}
\end{figure}

In Figure \ref{fig:error_distributions_odmd_dmdenkf} for Online DMD, the DMDEnKF and Hankel-DMDEnKF, we plot the distributions of 50 step ahead forecast mean run relative errors at both low and high levels of measurement noise. At low levels of noise, Online DMD's errors peak at a lower level than those of the DMDEnKF and Hankel-DMDEnKF, however as noise levels increase we see the peaks switch sides, and the DMDEnKF/Hankel-DMDEnKF become the better performing methods. At both noise levels, the peak in the DMDEnKF's error distribution is centred at the same value as the Hankel-DMDEnKF's peak, however it is less dense due to the additional probability mass stored in the long tail of the DMDEnKF's error distribution, which is not present in that of the Hankel-DMDEnKF.

These disproportionately large errors in the DMDEnKF distribution's tail occur when the spin-up DMD process fails to produce a model similar enough to the system's true dynamics. As briefly touched upon in the first synthetic example, if the spin-up DMD model is sufficiently inaccurate then it can stop the EnKF from effectively assimilating new data, leading to the catastrophic failure of the DMDEnKF. In this example, as the signal-to-noise ratio in the direction of the second-largest singular value was often low, an unfortunate random draw of the system dynamics ($\mathbf{A}$) and measurement noise $(\mathbf{v}_k)$ in the spin-up period could produce large errors in DMD's second eigenvalue. Empirically, using the interquartile range method to detect outlying forecast errors, this DMDEnKF failure occurred $5.5\%$ of the time for $\sigma =0.05$. The errors would persist in the filtering step as new data was processed, whereas other methods were able to recover from poor initial model estimates more effectively. The quality of the model produced by the initial DMD is dependent on the quality and volume of the spin-up data, hence this problem was exacerbated and occurred much more regularly at higher noise levels ($21.9\%$ of the time for $\sigma =0.5$). It could be mitigated somewhat by increasing the number of time steps in the spin-up stage as described at the end of Section 3.1, however similarly to Windowed TDMD as the system is assumed to be time varying there likely exists a point of negative returns once the spin-up period becomes too long due to the stationarity assumption of batch DMD becoming progressively more violated.

\begin{figure}[h]
     
     \begin{subfigure}[b]{0.48\textwidth}
         \centering
         \includegraphics[width=\textwidth]{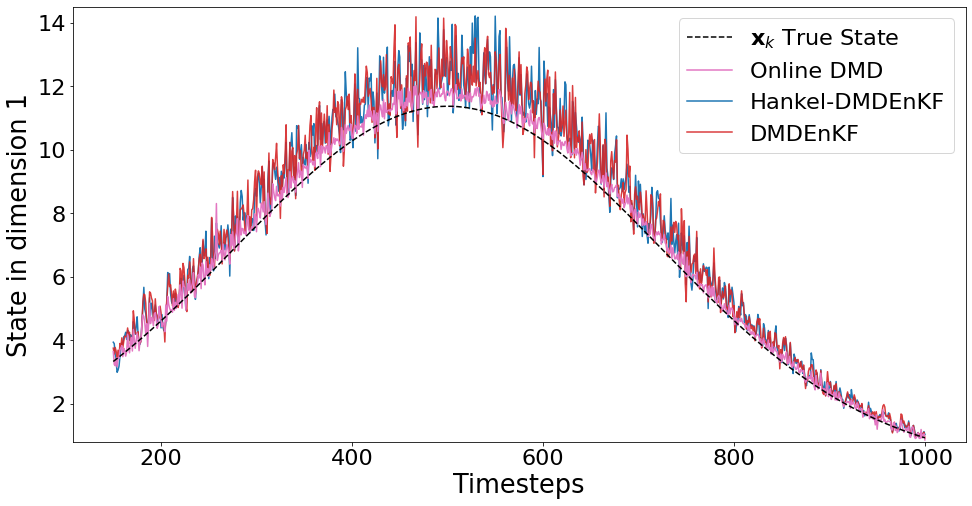} 
          \caption{\label{fig:0.05growth_decay_typical_trajectory}50 step ahead forecasts for $\sigma = 0.05$.}
          \end{subfigure}
          \hfill
    \begin{subfigure}[b]{0.48\textwidth}
         \centering
         \includegraphics[width=\textwidth]{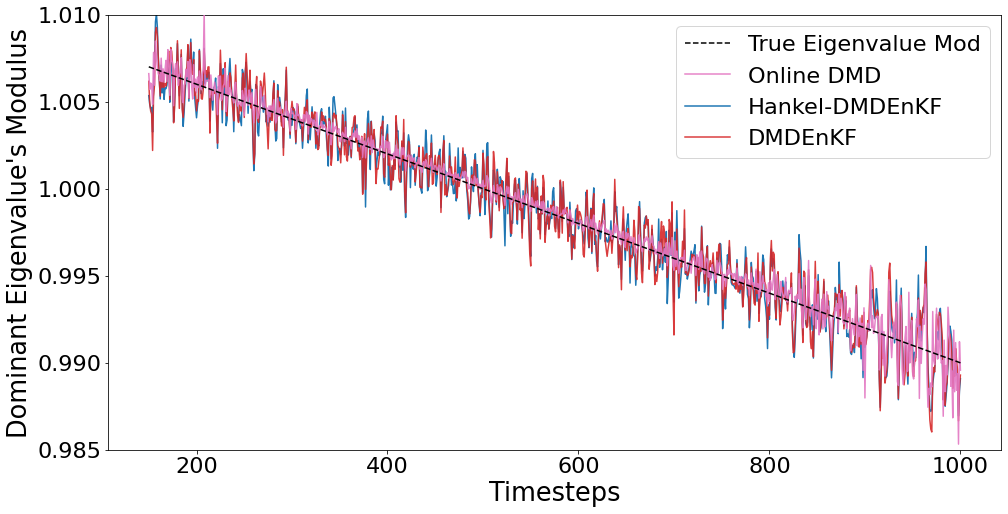} 
          \caption{\label{fig:0.05growth_decay_eigenvalue_tracking}Dominant eigenvalue's modulus for $\sigma = 0.05$.}
          \end{subfigure}
          \hfill
    \begin{subfigure}[b]{0.48\textwidth}
         \centering
         \includegraphics[width=\textwidth]{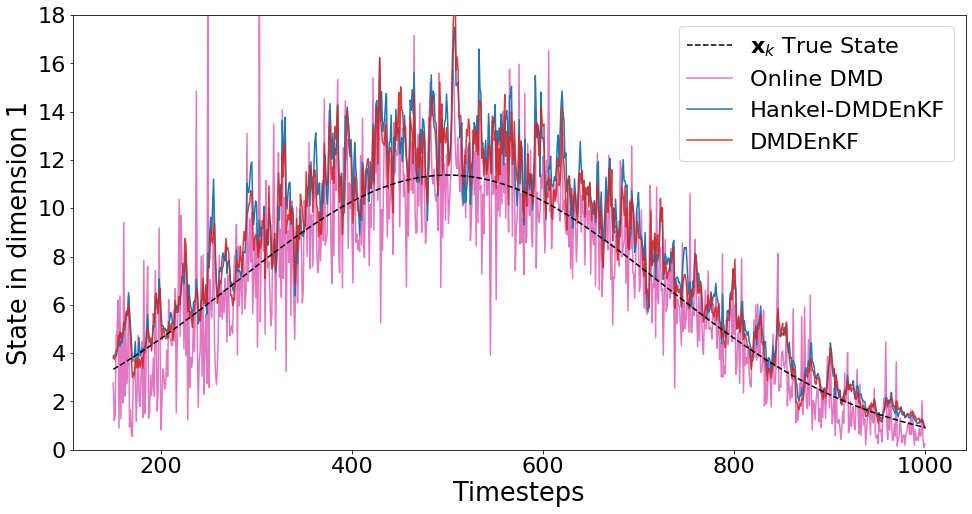} 
          \caption{\label{fig:0.5growth_decay_typical_trajectory}50 step ahead forecasts for $\sigma = 0.5$.}
          \end{subfigure}
          \hfill
         \begin{subfigure}[b]{0.48\textwidth}
         \centering
         \includegraphics[width=\textwidth]{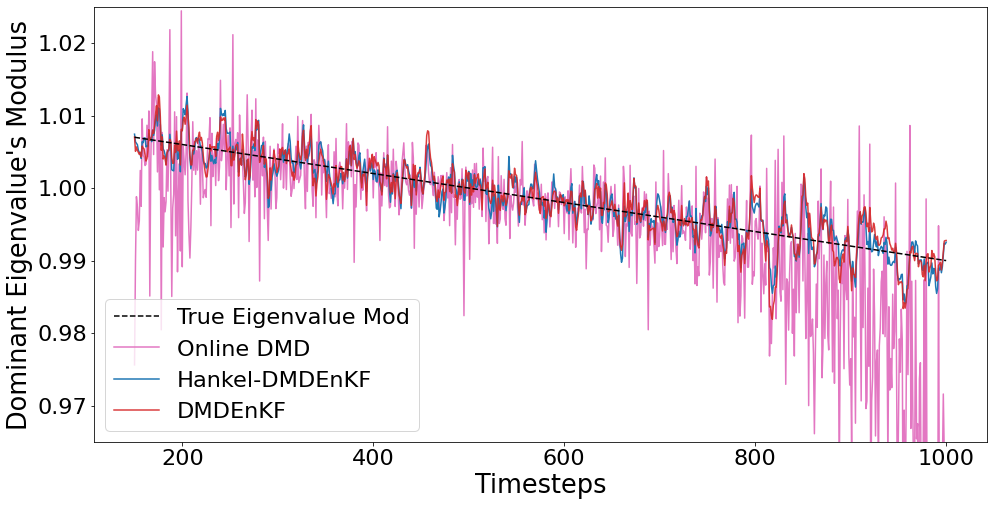}
          \caption{\label{fig:0.5growth_decay_eigenvalue_tracking}Dominant eigenvalue's modulus for $\sigma = 0.5$.}
          \end{subfigure}
          \caption{\centering\label{fig:growth_decay_trajectory_eigenvalue}Typical trajectories of the 50 step ahead forecasts for the value of the state's first dimension (left) and estimates of the dominant eigenvalue's current modulus (right) under low ($\sigma=0.05$) and high ($\sigma=0.5$) levels of measurement noise for Online DMD, the DMDEnKF and Hankel-DMDEnKF over the course of 1 run. Online DMD forecasts 50 steps ahead more accurately at low noise, and the DMDEnKF/Hankel-DMDEnKF more accurately at high noises, however when signal-to-noise ratio is low (at the start and end of the experiment) Online DMD's eigenvalue estimates become unstable.}
\end{figure}

Unlike the DMDEnKF, the Hankel-DMDEnKF and Online DMD do not suffer from a long tail in their error distributions, and perform consistently well over all 1000 runs. At both noise levels, their error distributions have a similar variance, with the Hankel-DMDEnKF's errors being slightly more tightly grouped than those of Online DMD. Hence, the average error is the main factor when differentiating between the method's performance in this example, meaning Online DMD is the preferred method at low noise and the Hankel-DMDEnKF (or DMDEnKF provided the spin-up DMD does not catastrophically fail) is more accurate at high noise. As both methods posses useful yet differing attributes, we generate a typical data trajectory (one for which the DMDEnKF does not fail) for both low and high measurement noise. We then investigate how each model's 50 step ahead forecasts and dominant eigenvalue estimates change over the course of each run, as shown in Figure \ref{fig:growth_decay_trajectory_eigenvalue}.

First, observing the low noise forecasts in Figure \ref{fig:0.05growth_decay_typical_trajectory}, it is clear Online DMD produces forecasts that are more robust and closer to the true state's value than those of the DMDEnKF and Hankel-DMDEnKF. This was to be expected, by virtue of Online DMD's lower average errors in the  error distributions of Figure \ref{fig:0.05_error_distributions_odmd_dmdenkf} at this noise level. As noise is increased, the forecasts in Figure \ref{fig:0.5growth_decay_typical_trajectory} show the DMDEnKF and Hankel-DMDEnKF becoming the more accurate methods, however Online DMD's forecasts remains fairly stable, and still appear to be a viable forecasting option.

Analysing the eigenvalue estimates in Figures \ref{fig:0.05growth_decay_eigenvalue_tracking} and \ref{fig:0.5growth_decay_eigenvalue_tracking}, we see that over the middle section of data where $k = (250,...,750)$, Online DMD is able to track the dominant eigenvalue effectively. However, at the beginning and end of the dataset when states and hence the signal component of each new data point is small relative to the measurement noise, Online DMD's eigenvalue estimates become progressively more unstable. In the low noise case this is not a problem, as Online DMD's estimates are significantly more accurate than those of the DMDEnKF/Hankel-DMDEnKF, so even in the poorly performing sections of the data it's estimates still better/match those of the DMDEnKF. For higher noise however, Online DMD provides significantly less robust estimates of the dominant eigenvalue at the start and end of the datasets than those generated by the DMDEnKF and Hankel-DMDEnKF. In the epidemiological context of an infectious disease outbreak, which this synthetic example attempts to mimic, scientists will often try to calculate the basic reproduction number ($R_0$) \cite{r0} using noisy data from the small number of initial infections. If $R_0 > 1$ the number of infections will grow exponentially if left unchecked, and if $R_0 < 1$ the number of infections will decay naturally to 0. Within this example, using the DMDEnKF/Hankel-DMDEnKF one can quickly determine that initially $R_0 > 1$ and take any required action thanks to the stability of it's early eigenvalue estimates, whereas it takes significantly longer and a higher level of infection for Online DMD to consistently determine if $R_0$ is above or below the growth/decay threshold.

\section{Seasonal Influenza-like Illness Application}
\subsection{Problem setup}
DMD based methods have previously been applied to infectious disease data \cite{epidemiology_use}. In this case, DMD modes can be viewed as stationary, spatial modes used to create a reduced order model in which only the amplitudes and frequencies are time varying \cite{epidemiology_rom_dmd}. Hence, modelling influenza-like illness (ILI) data is a prime potential application for the DMDEnKF/Hankel-DMDEnKF.

The CDC's ILINet data \cite{cdc_ili_definition} we will be using records the number of ILI General Practitioner (GP) consultations in the US each week, alongside the number of total GP consultations which can be used to normalize the ILI data. We use a subset of the data from the start of 2003, the first year when data is available all year round, to the end of 2018 as seen in Figure \ref{fig:ili_time_series}. We then split each week's data into demographics, consisting of 4 age groups (0-4, 5-24, 25-24, 65+) and 10 Health and Human Services (HHS) regions. Each region consists of the following locations:
\begin{itemize}

\item Region 1 - Connecticut, Maine, Massachusetts, New Hampshire, Rhode Island, and Vermont.

\item Region 2 - New Jersey, New York, Puerto Rico, and the U.S. Virgin Islands. 

\item Region 3 - Delaware, District of Columbia, Maryland, Pennsylvania, Virginia, and West Virginia.

\item Region 4 - Alabama, Florida, Georgia, Kentucky, Mississippi, North Carolina, South Carolina, and Tennessee.

\item Region 5 - Illinois, Indiana, Michigan, Minnesota, Ohio, and Wisconsin.

\item Region 6 - Arkansas, Louisiana, New Mexico, Oklahoma, and Texas.

\item Region 7 - Iowa, Kansas, Missouri, and Nebraska.

\item Region 8 - Colorado, Montana, North Dakota, South Dakota, Utah, and Wyoming.

\item Region 9 - Arizona, California, Hawaii, and Nevada.

\item Region 10 - Alaska, Idaho, Oregon, and Washington.

\end{itemize}

 Whilst ILI consultation data is available over all 40 of these strata, total GP consultation data is only provided by region. To generate an age breakdown for a region's total consultations we linearly interpolate using census data to approximate the US population's age demographics for a given week. We then allocate the region's total consultations to each age group based on the proportion of the total population they represent. This method assumes that all age groups have a similar likelihood of attending the GP's, which may be flawed but we believe it to be sufficient for the purpose of demonstrating the DMDEnKF on real-world data.

 \subsection{Building the spin-up DMD model}
\begin{figure}[!htbp]
\begin{subfigure}[b]{\textwidth}
         \includegraphics[width=\textwidth]{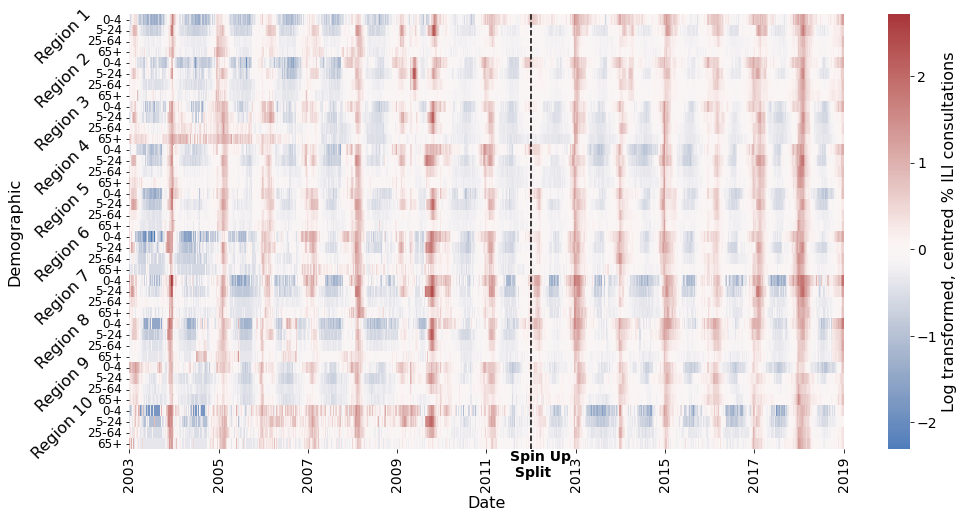}
     \end{subfigure}
     \caption{\centering\label{fig:ili_transformed_demo_data}ILI consultations as a percentage of total weekly GP consultations across 4 age brackets and the 10 HHS regions in the US, log transformed and centred. The peaks of varying size, timing and shape in Figure \ref{fig:ili_time_series} are visible here as vertical red areas of varying width and intensity that encompass most demographics.}
\end{figure}

 The format of the ILI data used in the DMDEnKF is thus a 40 dimensional vector for each week, recording ILI consultations as a percentage of total GP consultations over every demographic. This data exists in ${\mathbb R_{+}^{40}}$ however DMD computes modes in $\mathbb{R}^{40}$. Hence, to ensure the model's estimates are non-negative, we first transform the data by adding a small constant ($c=1$) then taking the natural logarithm of each element. For the Hankel-DMDEnKF, this transformed data is then delay-embedded with the previous 99 time steps ($d=100$) to form a state in $\mathbb{R}^{4000}$. We use data up to the end of 2012 as training data for the spin-up DMD processes of the DMDEnKF/Hankel-DMDEnKF detailed in Step 1 of the DMDEnKF algorithm, and then filter the remaining years from 2013-2018. The transformed, centred data with the split where the spin-up process ends and the filtering begins marked is shown in Figure \ref{fig:ili_transformed_demo_data}.

We initially choose to truncate to 8 DMD modes for the DMDEnKF to demonstrate the method. We discuss the effect of changing the truncation on the DMDEnKF method below, but at 8 DMD modes approximately the amount of additional variance in the data that is retained by keeping more modes diminishes significantly. This is evidenced by the ``elbow" seen in the cumulative variance plot of Figure \ref{fig:svd_variance_elbow}, at the point where the graph transitions from rapidly increasing in variance with $r$ to a more gradual ascent. We also truncate the Hankel-DMDEnKF to 8 DMD modes, to allow for a more direct comparison between the two variants.

\begin{figure}[!htbp]
\begin{subfigure}[b]{.42\textwidth}
         \includegraphics[width=\textwidth]{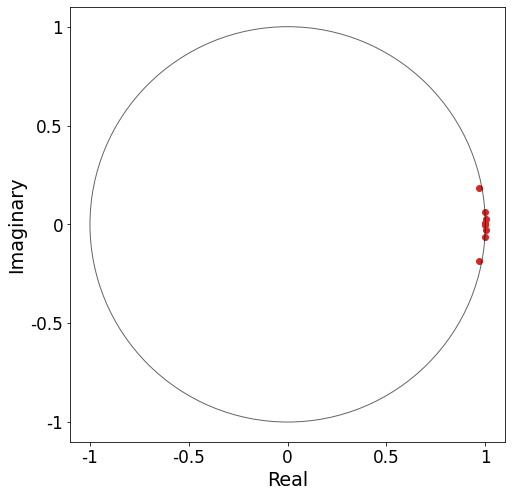}
              \caption{\centering\label{fig:ili_spec_and_modesa}DMD Eigenvalue Spectrum.}
     \end{subfigure}
\hfill
\begin{subfigure}[b]{.57\textwidth}
         \includegraphics[width=\textwidth]{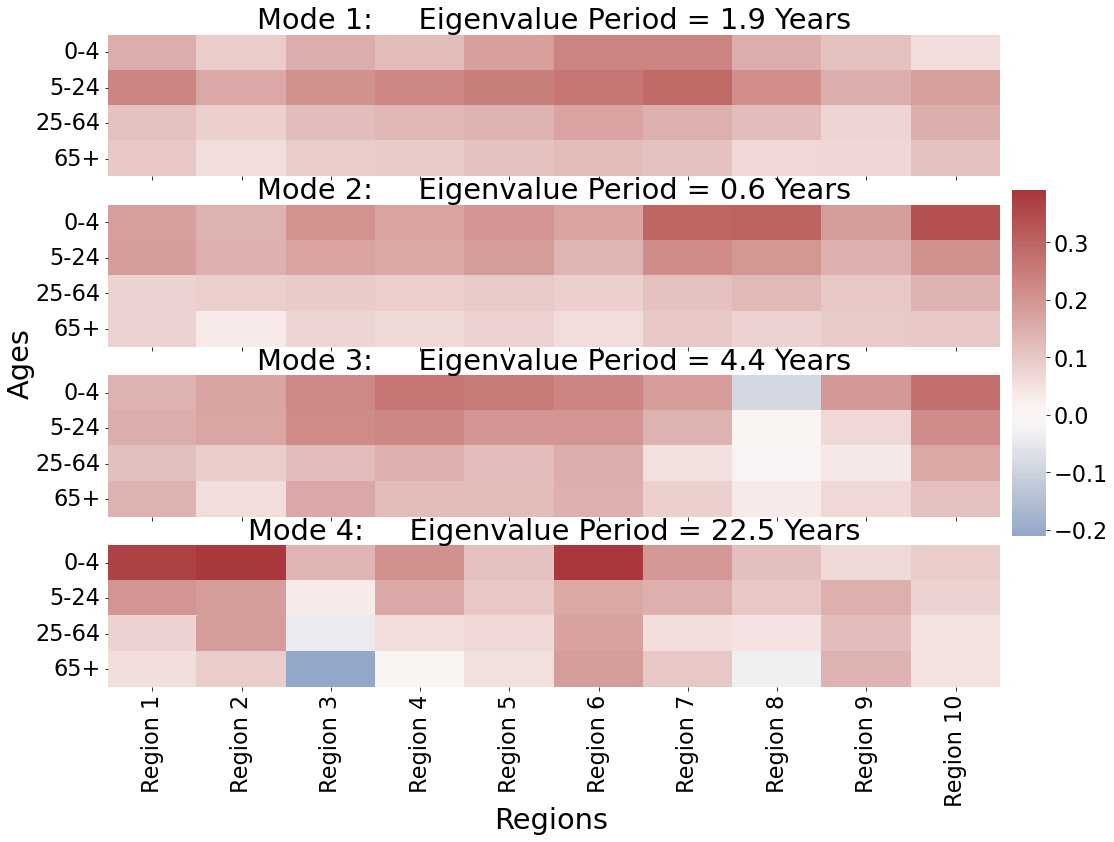}
              \caption{\centering\label{fig:ili_spec_and_modesb}Dominant DMD Modes.}
     \end{subfigure}
     
     \begin{subfigure}[b]{.42\textwidth}
         \includegraphics[width=\textwidth]{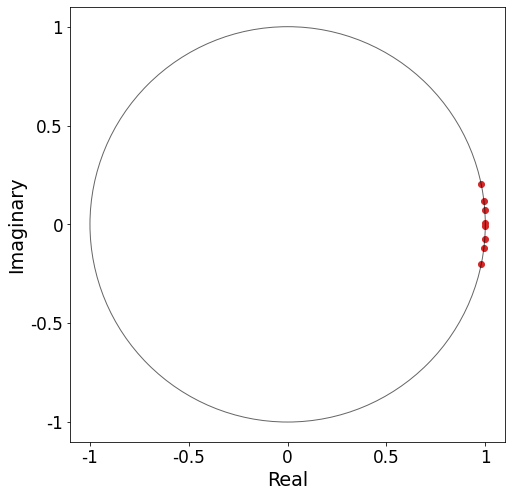}
              \caption{\centering\label{fig:ili_spec_and_modesc}Hankel-DMD Eigenvalue Spectrum.}
     \end{subfigure}
\hfill
\begin{subfigure}[b]{.57\textwidth}
         \includegraphics[width=\textwidth]{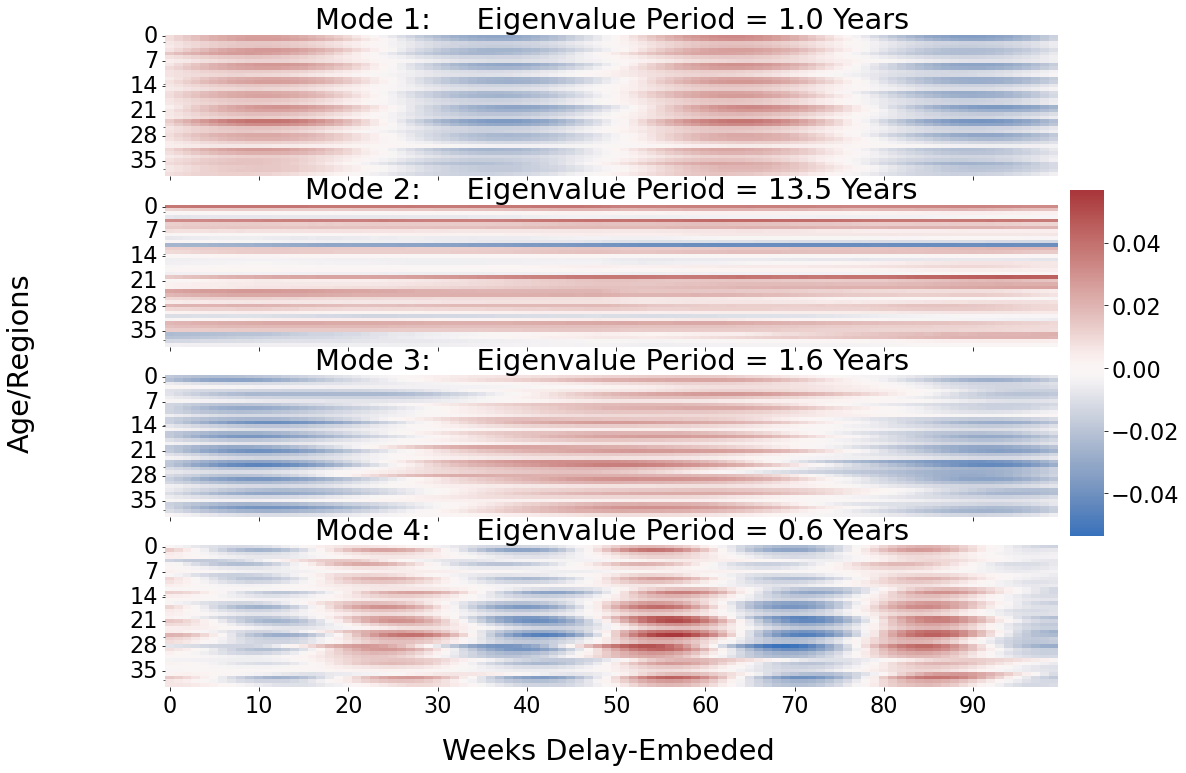}
              \caption{\centering\label{fig:ili_spec_and_modesd}Dominant Hankel-DMD Modes.}
     \end{subfigure}
     \caption{\centering\label{fig:ili_spec_and_modes}Eigenvalue Spectrum (left) and DMD modes in descending order of dominance (right) generated by the DMD (top)/Hankel-DMD (bottom) applied to the data in Figure \ref{fig:ili_transformed_demo_data} up to the spin-up date. In both cases, all eigenvalues lie approximately on the unit circle, and dominant modes feature the same sign for most demographics with a magnitude that varies with age. The DMD modes are more interpretable, but Hankel-DMD identifies the mode with period 1 year, which DMD does not.}
\end{figure}

The spectrum and dominant DMD/Hankel-DMD modes associated with each frequency identified by the spin-up DMD processes can be seen in Figure \ref{fig:ili_spec_and_modes}. All eigenvalues shown in Figures \ref{fig:ili_spec_and_modesa} and \ref{fig:ili_spec_and_modesc} had a modulus of $\sim 1$, meaning in both cases each mode was expected to persist in the data without growing exponentially. The major difference between the two methods spectra is that Hankel-DMD identifies the most dominant mode to have a period of one year, whereas DMD does not detect any modes with this period. Annual peaks in ILI consultations occurring at a relatively similar time each year indicates that the data contains a strong mode of period one year, and this is supported by Fourier analysis \cite{fourier_transform} which also identifies one year as the dominant period. Hence, DMD is missing the yearly mode present in the data which Hankel-DMD is able to detect, and this is likely due to Hankel-DMD's aforementioned enhanced robustness to measurement noise. There are two clear patterns in the structure of the dominant DMD and Hankel-DMD modes seen in Figures \ref{fig:ili_spec_and_modesb} and \ref{fig:ili_spec_and_modesd}. Firstly, their strata generally move together. This is shown by the vast majority of entries for each DMD mode, and entries within the same delay-embedded week (denoted by a vertical slice through the mode) for Hankel-DMD modes, sharing the same sign. This implies that the percentage of ILI consultations increases and decreases at a similar time across all ages and regions. Secondly, the variance is higher for the younger age groups. This is demonstrated by the absolute value of elements in the top two rows of each region generally being larger than those in the bottom two. From Figure \ref{fig:ili_spec_and_modesb}, this is visible trivially for the DMD modes. In Figure \ref{fig:ili_spec_and_modesd}, age groups are arranged in ascending order for each region, so this effect is evidenced in the Hankel-DMD modes by the presence of more intensely coloured horizontal lines of width 2, followed by less intensely coloured lines of width 2 repeating over each of the 10 regions. This indicates that there are sharper peaks and deeper troughs in the percentage of ILI consultations for the young, while the rates for those 25 and over remain more stable.

\subsection{Applying the filter}
\begin{figure}[!htbp]
\begin{subfigure}[b]{\textwidth}
         \includegraphics[width=\textwidth]{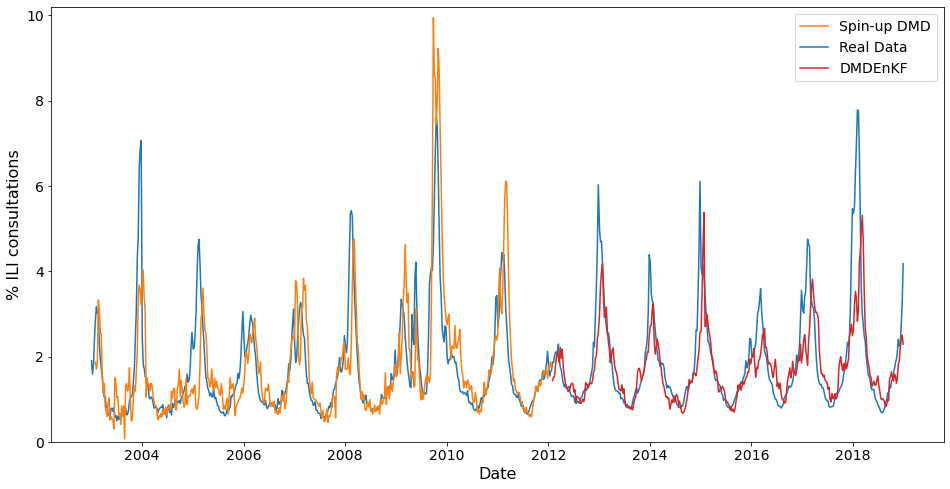}
         \caption{\centering\label{fig:ili_overall_forecasta}DMDEnKF 4-week ahead forecast.}
     \end{subfigure}
     
     \begin{subfigure}[b]{\textwidth}
         \includegraphics[width=\textwidth]{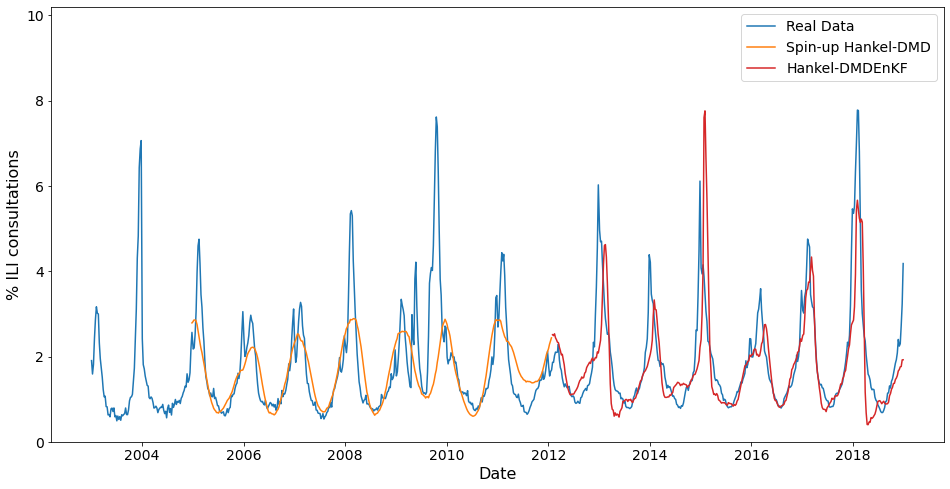}
             \caption{\centering\label{fig:ili_overall_forecastb}Hankel-DMDEnKF 4-week ahead forecast.}
     \end{subfigure}
     \caption{\centering\label{fig:ili_overall_forcast}ILI consultations as a percentage of total weekly GP consultations forecast 4 weeks ahead using the DMDEnKF (top) and Hankel-DMDEnKF (bottom). The DMD reconstruction captures the shape of the data well but is unstable, whereas the Hankel-DMD reconstruction is less sensitive to noise but over-smooths. The DMDEnKF and Hankel-DMDEnKF  forecasts help reduce these issues present in their respective reconstructions, but both suffer from some degree of lag in their predictions.}
\end{figure}

The filtering steps of the DMDEnKF/Hankel-DMDEnKF are then applied over the remaining data using the spatial and temporal modes from the spin-up DMD and spin-up Hankel-DMD respectively. Producing a 4-week ahead ILI forecast for the ILINet data that consistently outperforms a simple historical baseline prediction is difficult even for state-of-the-art models \cite{ilinet_model_comparison}. As such, to test the DMDEnKF/Hankel-DMDEnKF we use a forecast horizon of 4 weeks when making predictions. In Figure \ref{fig:ili_overall_forcast}, the DMDEnKF and Hankel-DMDEnKF's forecasting of total ILI consultations as a percentage of total weekly GP consultations can be seen in full. Up until 2012, we generate 4-week ahead reconstructions using the spin-up DMD models only, estimating each state by taking the state measurement from 4 weeks prior, then iteratively applying the model to it 4 times. The 4-week ahead DMD reconstruction in Figure \ref{fig:ili_overall_forecasta} captures more fluctuations in the data than that of Hankel-DMD, however these high frequency fluctuations can also indicate the effect of noise in the measurements. The Hankel-DMD reconstruction shown in Figure \ref{fig:ili_overall_forecastb} is much less sensitive to noise, although fails to identify the sharper peaks in the data, which suggest it may be over-smoothing. From 2012 onwards the filtering begins, and forecasts are generated as described in the DMDEnKF algorithm using equation \eqref{eqn:x=PLzPz}. The DMDEnKF forecasts become significantly more stable, while the Hankel-DMDEnKF forecasts improve in capturing the true shape of the data, however both suffered from some degree of lag in their predictions.

\begin{figure}[!htbp]
\begin{subfigure}[b]{\textwidth}
         \includegraphics[width=\textwidth]{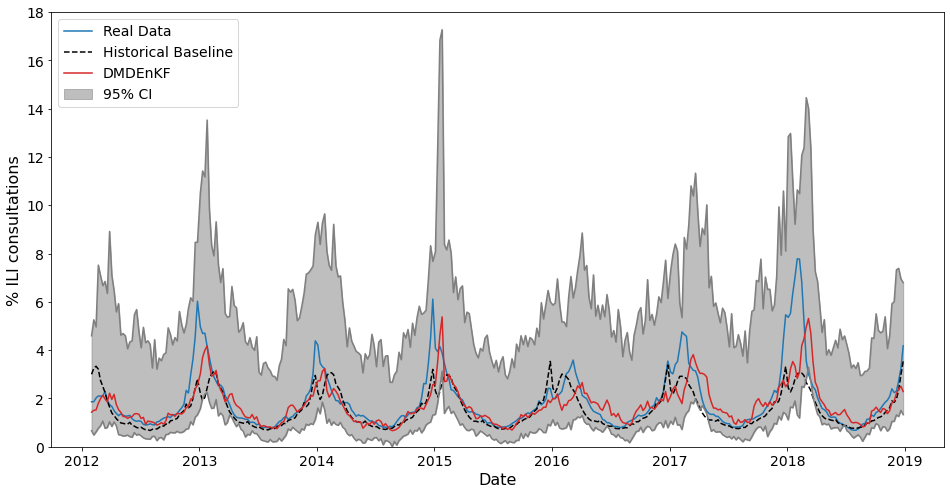}
     \end{subfigure}
     
 \begin{subfigure}[b]{\textwidth}
     \includegraphics[width=\textwidth]{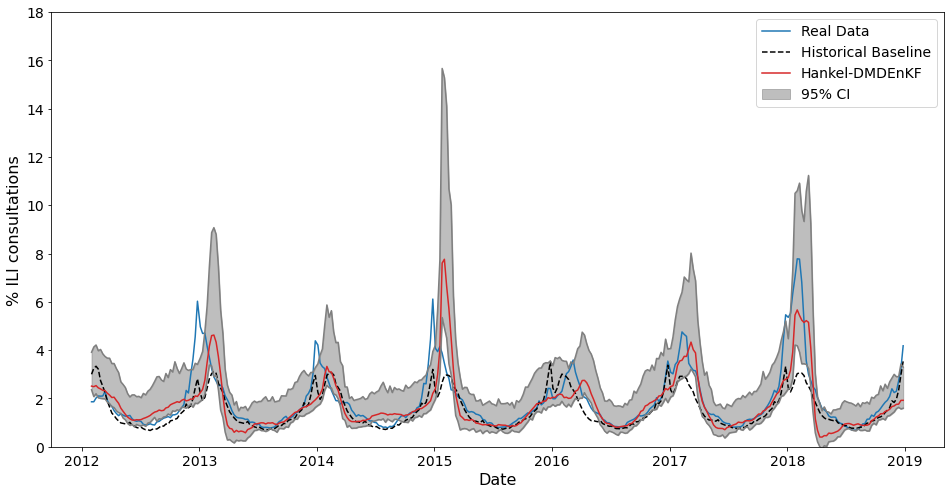}
 \end{subfigure}
     \caption{\centering\label{fig:ili_zoomed_forcast}ILI consultations as a percentage of total weekly GP consultations, forecast 4 weeks ahead using the DMDEnKF (top) and Hankel-DMDEnKF (bottom). A 95\% confidence interval for each forecast, and historical baseline predictions are also shown. The Hankel-DMDEnKF forecasts are smoother than those of the DMDEnKF, but both forecasts contain some lag. The real data always lies within the DMDEnKF's confidence interval but not the Hankel-DMDEnKF's, however this is likely due to the DMDEnKF's confidence interval being significantly wider than that of the Hankel-DMDEnKF.}
\end{figure}

During this second section of the data, the models are producing actual forecasts, as the DMDEnKFs only have access to data up to 4 weeks prior to the prediction target's date. Hence, it is in this section of the data we compare the models' performance against that of the historical baseline. The historical baseline prediction was created in a similar manner to that used in \cite{second_year_ilinet_comp}, taking ILI consultation rates from the same week of every previous year in the data (excluding the pandemic year of 2009) and then producing a probability distribution for the current week's consultations via Gaussian kernel density estimation (KDE) \cite{kernel_density_estimation}. KDE Bandwidths were determined using Silverman's rule of thumb \cite{silverman_density_estimation}, and when point estimates were required they were taken as the median of the distribution. The results of the comparisons can be seen in Figure \ref{fig:ili_zoomed_forcast} and Table \ref{table:dmdenkf_mse}. Here it is worth noting that although we use data dated 4 weeks prior to the prediction date, in reality this data is often subject to revisions so the ILINet data as it currently stands would not necessarily be available in real time \cite{ilinet_model_comparison}.

\subsection{Evaluating the DMDEnKF's performance}
Figure \ref{fig:ili_zoomed_forcast} demonstrates graphically the 4-step ahead DMDEnKF, Hankel-DMDEnKF and historical baseline forecasts. The DMDEnKFs more successfully capture the shape and height of each flu season's peak, however tend to predict the peaks late, whilst the historical baseline consistently underpredicts the peak rates but is fairly accurate on the timings. The Hankel-DMDEnKF's forecasts are smoother than those of the DMDEnKF, however do not capture smaller details within the shape of the peaks. We also plot the 95\% confidence intervals for the DMDEnKF and Hankel-DMDEnKF's forecasts in Figure \ref{fig:ili_zoomed_forcast}, generated using the ensemble that is maintained and propagated in the EnKF framework. At all times, the real data lies within the DMDEnKF’s
confidence interval, which is not true for the Hankel-DMDEnKF. The DMDEnKF's confidence interval is significantly wider than that of the Hankel-DMDEnKF, and this is due to Hankel-DMD's robustness to noise, meaning that when the ensemble is propagated through the model, a large amount of probability mass is concentrated in a small area of the state space. This then leads to the Hankel-DMDEnKF underestimating the uncertainty in the system, and hence some real data values falling outside the boundaries of it's 95\% confidence interval.

\begin{table}[h]
\begin{center}
 \begin{tabular}{||c c c||} 
 \hline
 Forecast & Log Score & Mean Squared Error\\ [0.5ex] 
 \hline\hline
 Historical Baseline & 0.28 & 1.24 \\
 \hline
 1-week ahead DMDEnKF & 0.49 & 0.33 \\ 
 \hline
 2-week ahead DMDEnKF & 0.38 & 0.61 \\
 \hline
 3-week ahead DMDEnKF  & 0.32 & 0.87 \\
 \hline
 4-week ahead DMDEnKF  & 0.27 & 1.16 \\
 \hline
 1-week ahead Hankel-DMDEnKF & 0.41 & 0.49 \\ 
 \hline
 2-week ahead Hankel-DMDEnKF & 0.33 & 0.70 \\
 \hline
 3-week ahead Hankel-DMDEnKF & 0.29 & 0.97 \\
 \hline
 4-week ahead Hankel-DMDEnKF & 0.23 & 1.26 \\
 \hline
\end{tabular}
\end{center}
\caption{\centering The log scores and mean squared errors for the DMDEnKF and Hankel-DMDEnKF with differing forecast horizons, and the historical baseline prediction. The DMDEnKF achieves a higher log score and mean squared error than the historical baseline for forecast horizons up to 4 weeks ahead, where it attains a similar level of forecast skill. The Hankel-DMDEnKF consistently underperforms against the DMDEnKF in both metrics over these short forecast horizons. Scores are calculated over the 6 flu seasons from 2012/13 to 2017/18.}
\label{table:dmdenkf_mse}
\end{table}

To numerically compare performance, we used metrics designed for the Forecast the Influenza Season Collaborative Challenge (FluSight), in which multiple teams would submit predictions about the weekly ILINet consultation rates for the upcoming flu season at a national and HHS regional level \cite{first_year_ilinet_comp}, \cite{second_year_ilinet_comp}. The FluSight challenge evaluated models abilities to generate 1-4-week ahead forecasts known as short-term targets over the course of a flu season, as well as other longer term targets, known as seasonal targets, before the season had begun. The DMDEnKF is primarily intended to be a tool for tracking and short-term forecasting, hence we focus on forecasting these short-term targets only. For this purpose we used two different metrics, the log probability measure (log score) slightly adjusted from the FluSight challenge as used in \cite{ilinet_model_comparison} and the mean squared error due to its popular use in regression problems. The log score represents the geometric average probability of each model's prediction being accurate, with accuracy deemed as a forecast within $+/-0.5$ of the true ILI consultation rate. The higher the log score, the better the forecast. Metrics are calculated from week 40 to week 20 of the following year to prioritize evaluation of forecasts during the flu season, and we use the 6 full seasons from 2012/13 to 2017/18.

The results for the historical baseline prediction and DMDEnKF/Hankel-DMDEnKF's forecasts at a national level can be seen in Table \ref{table:dmdenkf_mse}. As one would expect, the accuracy of both DMDEnKFs degrade as they make predictions further into the future. The DMDEnKF achieves a higher log score and mean squared error than the historical baseline for forecast horizons up to 4 weeks ahead, where it attains a similar level of forecast skill. For forecasts of 5 or more weeks ahead, the DMDEnKF is unable to outperform the historical baseline in either metric. The Hankel-DMDEnKF consistently underperforms against the DMDEnKF in both metrics over these short forecast horizons. The top 3 statistical models and top 3 mechanistic models in the FluSight challenge achieved log scores of 0.32 and 0.3 respectively for their 4-week ahead forecasts, hence the DMDEnKF has lower (but comparable) forecasting skill than current state of the art ILI models. As the forecast horizon is extended up to 12 weeks ahead, the DMDEnKF's forecast scores continue to decrease monotonically, whereas the Hankel-DMDEnKF's log scores for 9-12 weeks ahead are no worse than those for 5-8 weeks ahead. As such, the DMDEnKF is preferred for short-term forecasting, while the Hankel-DMDEnKF is considered superior when forecasting over longer timescales.

 \begin{figure}[!htbp]

    \begin{subfigure}[b]{0.48\textwidth}
         \centering
         \includegraphics[width=\textwidth]{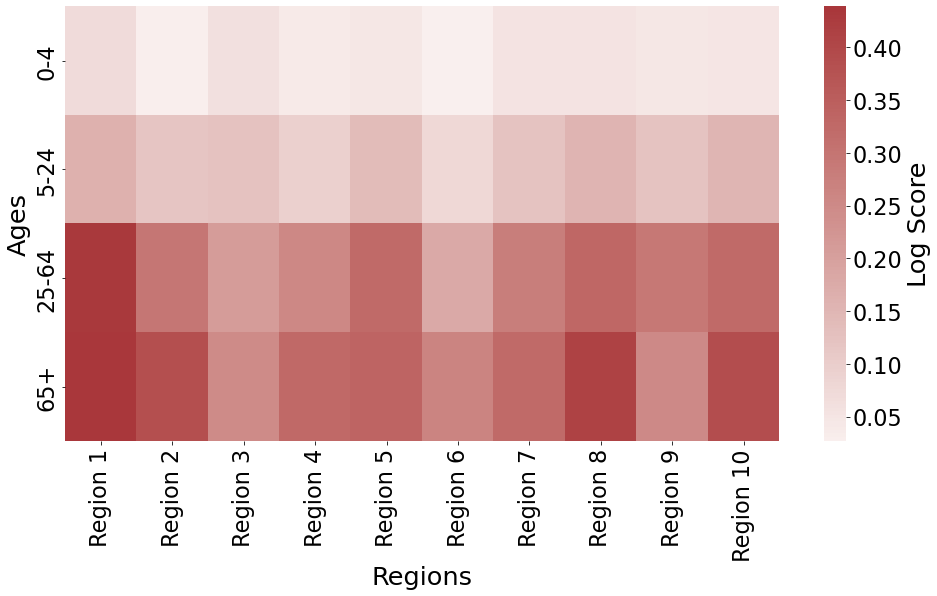} 
          \caption{\label{fig:age_region_log_scoresa}\centering DMDEnKF 4-week ahead forecast. \protect\newline }
          \end{subfigure}
          \hfill
         \begin{subfigure}[b]{0.48\textwidth}
         \centering
         \includegraphics[width=\textwidth]{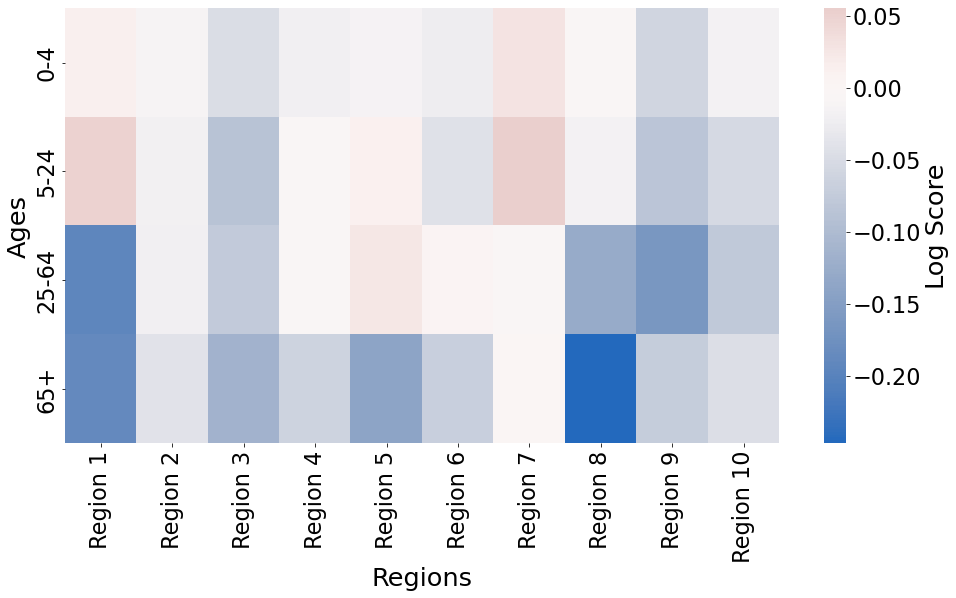}
          \caption{\label{fig:age_region_log_scoresb}\centering DMDEnKF 4-week ahead forecast - historical baseline.}
          \end{subfigure}
          \caption{\centering\label{fig:age_region_log_scores}Log scores over all ages and regions for the DMDEnKF's 4-week ahead forecast (left), followed by those same scores with the log scores of the historical baseline prediction subtracted (right). The Hankel-DMDEnKF scored similarly to the DMDEnKF across all ages and regions, so we do not include its breakdown to avoid redundancy. In the top figure, the generally increasing intensity of red as one moves down the age groups shows the DMDEnKF performing more accurately for older age groups. The bottom figure's varying areas of red/blue shows the DMDEnKF/historical baseline vary in superiority of forecasting skill depending on the age and region being forecast, with the historical baseline scoring more highly for most regions.}
\end{figure}

Figure \ref{fig:age_region_log_scores} shows the log scores for the 4-week ahead DMDEnKF forecast, and how these compare to the scores attained by the historical baseline prediction at an age and regional level. The Hankel-DMDEnKF scored similarly to the DMDEnKF across all ages and regions, so its breakdown is rather similar to that of the DMDEnKF, and we do not include it to avoid redundancy. In the DMDEnKF's log scores, we see a major pattern in the older age groups scoring higher and hence being better predicted than the younger demographics. This pattern does not persist when the historical baseline's scores are removed, indicating it is a more general trait of the data as opposed to a specific quality of the DMDEnKF's modelling technique. There is also a significant difference in the predictability from region to region. For example, region 1 was the most predictable region for both the DMDEnKF and historical baseline, which is consistent with the findings in \cite{ilinet_model_comparison}. However, the DMDEnKF improved on the historical baseline's forecast for only two of the four age groups in this region. In \cite{ilinet_model_comparison} it was found that the most overall improvement gained by forecasting for a region using a model as opposed to the historical baseline prediction also occurred in region 1, so one would expect to see improvements by the DMDEnKF over the historical baseline in all four age groups. As log score is heavily influenced by the amount of uncertainty in a forecast, it is possible that the covariance matrices used in the DMDEnKF were too large for this region. Hence, setting the magnitude of the DMDEnKF's variance on a region by region basis could lead to better results and more accurate forecasts. Region 6 was the worst forecast region by the DMDEnKF, and the historical baseline predictions were slightly more accurate. Again, this is consistent with \cite{ilinet_model_comparison} where region 6 was the lowest scoring region for the models. In that work however, region 6 experienced the second most improvement by using a model over the historical baseline prediction. Hence, for this region a larger variance within the DMDEnKF may have been more appropriate to account for its extra unpredictability, further supporting the idea of varying the variance by region in the future.

\subsection{Varying the truncation rank}
\begin{figure}[!htbp]
\begin{subfigure}[b]{.48\textwidth}
         \includegraphics[width=\textwidth]{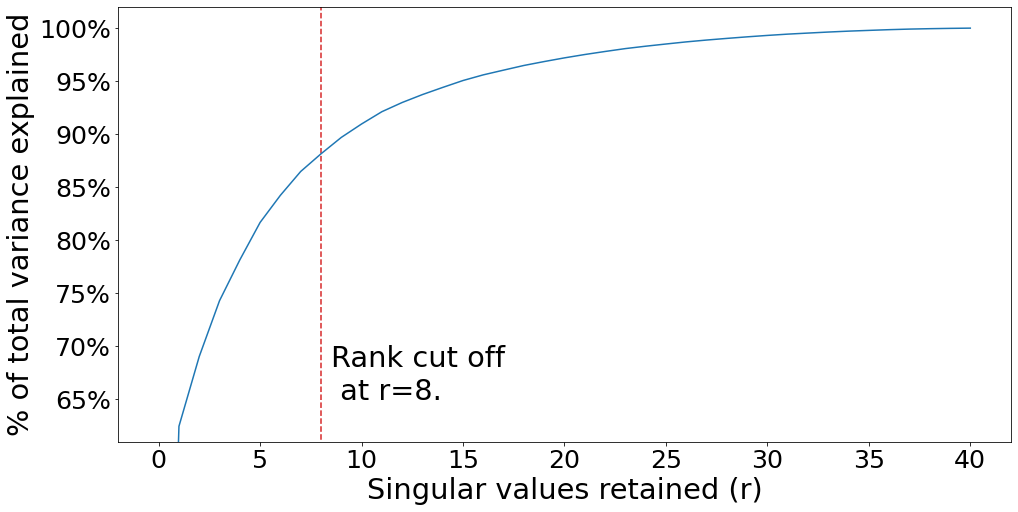}
                \caption{\centering\label{fig:svd_variance_elbow}$\%$ of total variance in the data.\protect\newline}
\end{subfigure}
\hfill     
\begin{subfigure}[b]{.48\textwidth}
         \includegraphics[width=\textwidth]{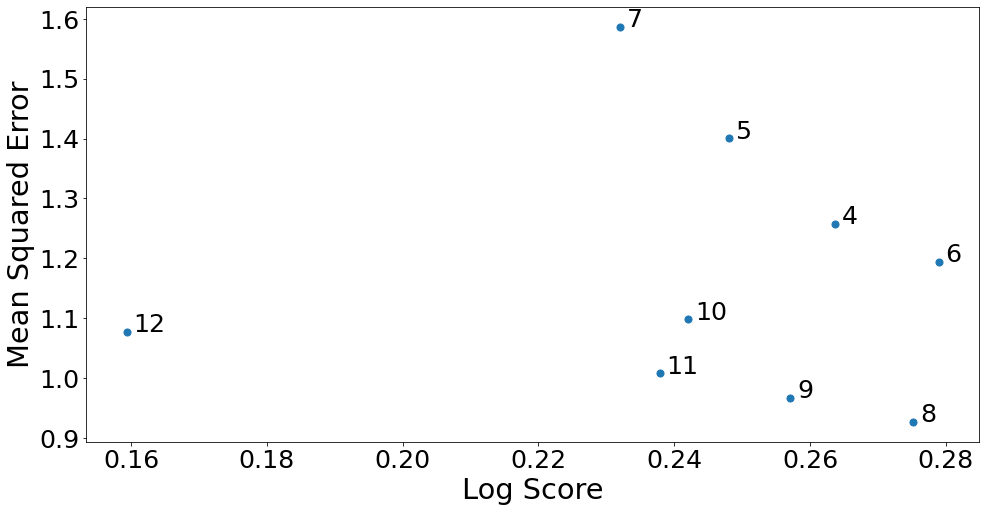}
        \caption{\centering\label{fig:rank_errors}DMDEnKF log score/mean squared errors for $r=4,..., 12$.}
     \end{subfigure}
     
     \begin{subfigure}[b]{.48\textwidth}
         \includegraphics[width=\textwidth]{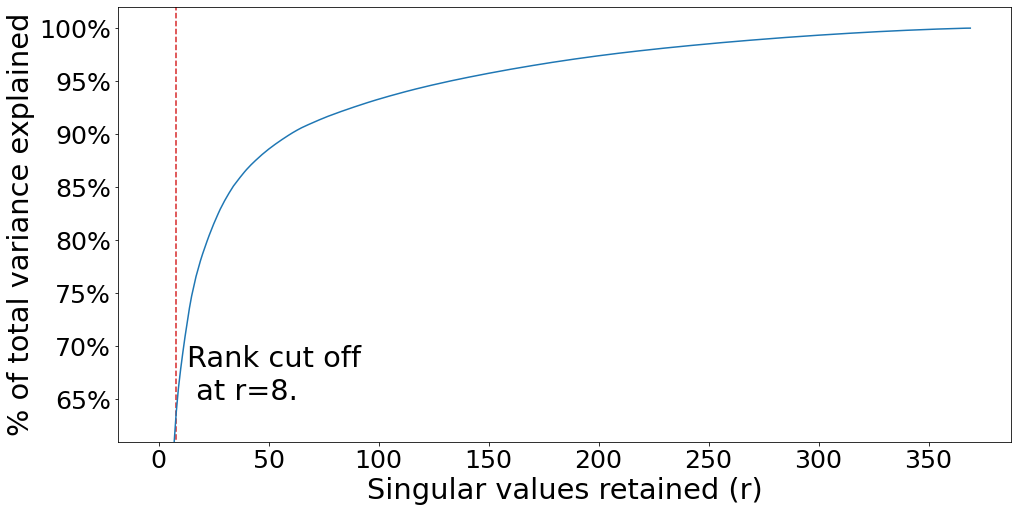}
                \caption{\centering\label{fig:hankel_svd_variance_elbow}$\%$ of total variance in the delay-embedded data.\protect\newline}
\end{subfigure}
\hfill     
\begin{subfigure}[b]{.48\textwidth}
         \includegraphics[width=\textwidth]{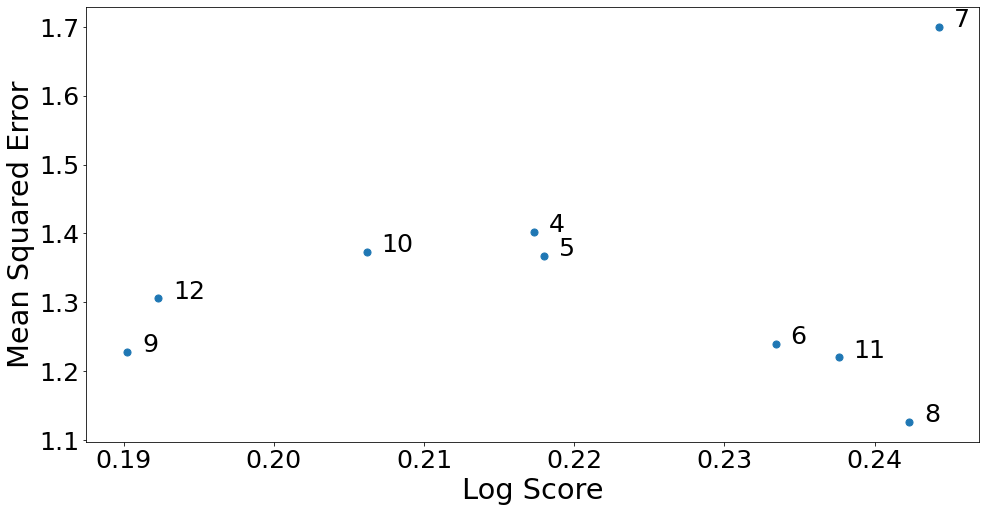}
        \caption{\centering\label{fig:hankel_rank_errors}Hankel-DMDEnKF log score/mean squared errors for $r=4,..., 12$.}
     \end{subfigure}
     \caption{\centering\label{fig:svd_var_elbow_rank_errors}On the left, the $\%$ of the total variance in the data (top) and delay-embedded data (bottom), dependent on the number of singular values that are retained ($r$). An ``elbow" in the data occurs around $r=8$ where we choose to truncate, however determining the exact position of the ``elbow" is subjective and could be considered anywhere from $r=4$ to $r=12$. On the right, the log score and mean squared errors for 4-step ahead forecasts generated using the DMDEnKF (top) and Hankel-DMDEnKF (bottom) with differing values of $r$. In both cases, log score is maximised and mean squared error minimised for $r = 8$.}
\end{figure}

Having analysed the DMDEnKF and Hankel-DMDEnKF's ILI forecasting with 8 DMD modes, we now investigate the effect different truncation ranks ($r$) have on their performance in Figure \ref{fig:svd_var_elbow_rank_errors}. From Figure \ref{fig:svd_variance_elbow}, the subjective process of identifying an ``elbow" in the data could lead an observer to determine a suitable rank for truncation as low as 4 or as high as 12. Application of the algorithm of Gavish and Donoho for identifying the optimal truncation threshold \cite{optimal_trunc} also finds the truncation rank to be 12, hence we will focus on investigating values of $r$ in the interval from 4 to 12.

Figures \ref{fig:rank_errors} and \ref{fig:hankel_rank_errors} show how the metrics we use to measure the DMDEnKF and Hankel-DMDEnKF's forecasting skill vary with $r$. An ideal forecast will have a high log score reflecting a relatively tight and accurate probability distribution, with a low mean squared error indicating a point estimate close to the true percentage of ILI consultations. For both methods, log score is maximised and mean squared error minimised by $r = 8$, indicating this is the optimal rank to truncate at for our models. For $r=4$, we have the simplest model tested, hence it has a low degree of uncertainty resulting in a relatively high log score, however is too simple to properly model the system so receives a high mean squared error. By increasing $r$, we allow the DMDEnKFs more freedom to capture complexity within the system, resulting in a more accurate representation of the true dynamics and hence a generally lower mean squared error. When the number of eigenvalues is increased too far however, it begins modelling elements of the noise in the system which negatively impacts future predictions, as seen in the increase in mean squared errors for $r > 8$. The additional freedom afforded to the DMDEnKF by increasing $r$ also means the model contains more parameters, each of which have an associated degree of uncertainty. This causes the overall forecast's probability distribution to become more spread out, and when no longer offset by the increased model accuracy up to $r = 8$, reduces the forecasts log score.

\section{Conclusion}

To conclude, we have defined two new algorithms, the DMDEnKF and Hankel-DMDEnKF, that combine dynamic mode decomposition and Hankel dynamic mode decomposition respectively with ensemble Kalman filtering, to update state and temporal mode estimates of a dynamical system as new data becomes available. When applied to simple, synthetic systems with a time varying parameter and low measurement noise, the DMDEnKFs performed similarly to other iterative DMD variants tested in tracking the system's time varying parameter and forecasting future states. As measurement noise was increased, the DMDEnKFs outperformed the other methods tested in both metrics, and the Hankel-DMDEnKF produced more stable forecasts than those of the DMDEnKF. Both DMDEnKFs achieved similar performance levels to their equivalent DMD Particle Filters (an alteration to the DMDEnKF algorithms where the ensemble Kalman filters were switched for Particle Filters), while requiring significantly fewer ensemble members. When forecasting influenza-like illness across age groups and HHS regions in the US using data from the CDC, the DMDEnKF produced more accurate forecasts than a historical baseline prediction up to 3 weeks ahead, and forecasts approximately as accurate 4 weeks ahead. The Hankel-DMDEnKF produced less accurate forecasts for these short-term targets than the DMDEnKF, however in general it's forecasts were more stable. Also, the Hankel-DMDEnKF was able to identify the presence of a mode with period 1 year, which is strongly visible in the data, yet not identified by the DMDEnKF. Both DMDEnKFs exhibited lower forecasting skill than current state of the art influenza-like illness models.

A natural extension of the DMDEnKF would be to apply extended/kernel DMD in the spin-up DMD phase, allowing the algorithm to be used more effectively on dynamical systems that act nonlinearly in their measured states. Instead of taking the observed values alone as the system's state $\mathbf{x}_k$, these variants use for the state a collection of functions on the observables $\mathbf{g}(\mathbf{x}_k)$, which often increases the state's dimension $n$. The EnKF is well suited to this pairing, as it scales more computationally efficiently in the state dimension than other Kalman filtering methods \cite{enkf}. The best choice of the collection of functions $\mathbf{g}(\mathbf{x}_k)$ as an embedding for nonlinear systems so that DMD may be effectively utilized is an interesting area of future work. Many methods have been developed that propose ways of generating $\mathbf{g}(\mathbf{x}_k)$, for example using deep learning \cite{dmd_autoencoder} or reservoir computing \cite{dmd_reservoir_computing}, and this remains a promising avenue for future work.

\section*{Code availability}
Codes used to produce the results in this paper are available at https://github.com/falconical/DMDEnKF.

\section*{Data availability statement}
All data used to produce the results in this paper will be made available upon reasonable request.

\section*{Acknowledgements}
This work was supported by the UKRI, whose Doctoral Training Partnership Studentship helped fund Stephen Falconers PhD. He would also like to thank for their valuable discussions, Nadia Smith and Spencer Thomas from the National Physics Laboratory.

%\clearpage

\bibliographystyle{siam}
\bibliography{mybibfile}

\begin{thebibliography}{10}

\bibitem{eigeninfo}
{\sc H.~Abdi}, {\em The eigen-decomposition: Eigenvalues and eigenvectors},
  Encyclopedia of measurement and statistics,  (2007).

\bibitem{hankel-dmd}
{\sc H.~Arbabi and I.~Mezić}, {\em Ergodic theory, dynamic mode decomposition,
  and computation of spectral properties of the koopman operator}, SIAM Journal
  on Applied Dynamical Systems, 16 (2017), pp.~2096--2126.

\bibitem{kf_applications}
{\sc F.~Auger, M.~Hilairet, J.~Guerrero, E.~Monmasson, T.~Orlowska-Kowalska,
  and S.~Katsura}, {\em Industrial applications of the kalman filter: A
  review}, IEEE Transactions on Industrial Electronics, 60 (2013), p.~5458.

\bibitem{mechanistic_vsml}
{\sc R.~Baker, J.-M. Peña, J.~Jayamohan, and A.~Jérusalem}, {\em Mechanistic
  models versus machine learning, a fight worth fighting for the biological
  community?}, Biology Letters, 14 (2018), p.~20170660.

\bibitem{ili_mechanistic1}
{\sc D.~Balcan, V.~Colizza, B.~Gon{\c c}alves, H.~Hu, J.~J. Ramasco, and
  A.~Vespignani}, {\em Multiscale mobility networks and the spatial spreading
  of infectious diseases}, Proceedings of the National Academy of Sciences, 106
  (2009), pp.~21484--21489.

\bibitem{moore-pen}
{\sc J.~C.~A. Barata and M.~S. Hussein}, {\em The moore–penrose
  pseudoinverse: A tutorial review of the theory}, Brazilian Journal of
  Physics, 42 (2011), p.~146–165.

\bibitem{first_year_ilinet_comp}
{\sc M.~Biggerstaff, D.~Alper, M.~Dredze, S.~Fox, I.~C.-H. Fung, K.~S.
  Hickmann, B.~Lewis, R.~Rosenfeld, J.~Shaman, M.-H. Tsou, et~al.}, {\em
  Results from the centers for disease control and prevention’s predict the
  2013--2014 influenza season challenge}, BMC infectious diseases, 16 (2016),
  pp.~1--10.

\bibitem{second_year_ilinet_comp}
{\sc M.~Biggerstaff, M.~Johansson, D.~Alper, L.~C. Brooks, P.~Chakraborty,
  D.~C. Farrow, S.~Hyun, S.~Kandula, C.~McGowan, N.~Ramakrishnan, R.~Rosenfeld,
  J.~Shaman, R.~Tibshirani, R.~J. Tibshirani, A.~Vespignani, W.~Yang, Q.~Zhang,
  and C.~Reed}, {\em Results from the second year of a collaborative effort to
  forecast influenza seasons in the united states}, Epidemics, 24 (2018),
  pp.~26--33.

\bibitem{epidemiology_rom_dmd}
{\sc D.~A. Bistrian, G.~Dimitriu, and I.~M. Navon}, {\em Processing
  epidemiological data using dynamic mode decomposition method}, AIP Conference
  Proceedings, 2164 (2019), p.~080002.

\bibitem{fourier_transform}
{\sc R.~N. Bracewell and R.~N. Bracewell}, {\em The Fourier transform and its
  applications}, vol.~31999, McGraw-Hill New York, 1986.

\bibitem{ili_statistical2}
{\sc L.~Brooks, D.~Farrow, S.~Hyun, R.~Tibshirani, and R.~Rosenfeld}, {\em
  Flexible modeling of epidemics with an empirical bayes framework}, PLOS
  Computational Biology, 11 (2014).

\bibitem{neuroscience_use}
{\sc B.~W. Brunton, L.~A. Johnson, J.~G. Ojemann, and J.~N. Kutz}, {\em
  Extracting spatial–temporal coherent patterns in large-scale neural
  recordings using dynamic mode decomposition}, Journal of Neuroscience
  Methods, 258 (2016), p.~1–15.

\bibitem{opt_dmd}
{\sc K.~Chen, J.~Tu, and C.~Rowley}, {\em Variants of dynamic mode
  decomposition: Boundary condition, koopman, and fourier analyses}, Journal of
  Nonlinear Science, 22 (2012).

\bibitem{ili_modelling_review2}
{\sc J.-P. Chretien, D.~George, J.~Shaman, R.~Chitale, and F.~McKenzie}, {\em
  Influenza forecasting in human populations: A scoping review}, PloS one, 9
  (2014), p.~e94130.

\bibitem{ncdmd}
{\sc S.~T.~M. Dawson, M.~S. Hemati, M.~O. Williams, and C.~W. Rowley}, {\em
  Characterizing and correcting for the effect of sensor noise in the dynamic
  mode decomposition}, Experiments in Fluids, 57 (2016).

\bibitem{time_delay_pca}
{\sc A.~{de Cheveigné} and J.~Z. Simon}, {\em Denoising based on time-shift
  pca}, Journal of Neuroscience Methods, 165 (2007), pp.~297--305.

\bibitem{particle_filter_convergence}
{\sc P.~{Del Moral}}, {\em Nonlinear filtering: Interacting particle
  resolution}, Comptes Rendus de l'Académie des Sciences - Series I -
  Mathematics, 325 (1997), pp.~653--658.

\bibitem{da_use_medicine}
{\sc M.~D'Elia, L.~Mirabella, T.~Passerini, M.~Perego, M.~Piccinelli,
  C.~Vergara, and A.~Veneziani}, {\em Applications of variational data
  assimilation in computational hemodynamics}, Modeling, Simulation and
  Applications, 5 (2012).

\bibitem{pydmd}
{\sc N.~Demo, M.~Tezzele, and G.~Rozza}, {\em {PyDMD: Python Dynamic Mode
  Decomposition}}, The Journal of Open Source Software, 3 (2018), p.~530.

\bibitem{particle_filter_resampling_schemes}
{\sc R.~Douc and O.~Capp{\'e}}, {\em Comparison of resampling schemes for
  particle filtering}, in ISPA 2005. Proceedings of the 4th International
  Symposium on Image and Signal Processing and Analysis, 2005., IEEE, 2005,
  pp.~64--69.

\bibitem{particle_filter_0.5_effective_sample_size}
{\sc A.~Doucet and A.~Johansen}, {\em A tutorial on particle filtering and
  smoothing: Fifteen years later}, Handbook of Nonlinear Filtering, 12 (2009).

\bibitem{eckart-young}
{\sc C.~Eckart and G.~Young}, {\em The approximation of one matrix by another
  of lower rank}, Psychometrika, 1 (1936), pp.~211--218.

\bibitem{enkf}
{\sc G.~Evensen}, {\em The ensemble kalman filter: Theoretical formulation and
  practical implementation}, Ocean dynamics, 53 (2003), pp.~343--367.

\bibitem{ili_deaths_updated}
{\sc S.~I. Flu}, {\em Estimated influenza illnesses, medical visits,
  hospitalizations, and deaths averted by vaccination in the united states},
  Prevent,  (2007), pp.~2006--2007.

\bibitem{cdc_ili_definition}
{\sc C.~for Disease~Control and Prevention}, {\em Cdc u.s. influenza
  surveillance system: Purpose and methods}.
\newblock \url{https://www.cdc.gov/flu/weekly/overview.htm}, 2021.
\newblock Accessed: 2021-04-27.

\bibitem{optimal_trunc}
{\sc M.~Gavish and D.~L. Donoho}, {\em The optimal hard threshold for singular
  values is 4/sqrt(3)}, 2014.

\bibitem{particle_filter_og}
{\sc N.~Gordon, D.~Salmond, and A.~F.~M. Smith}, {\em Novel approach to
  nonlinear/non-gaussian bayesian state estimation}, IEE Proceedings F (Radar
  and Signal Processing), 140 (1993), pp.~107--113(6).

\bibitem{dmd_reservoir_computing}
{\sc M.~Gulina and A.~Mauroy}, {\em Two methods to approximate the koopman
  operator with a reservoir computer}, Chaos: An Interdisciplinary Journal of
  Nonlinear Science, 31 (2021), p.~023116.

\bibitem{odmd_github}
{\sc C.~W.~R. Hao~Zhang}, {\em Online dmd github repository}.
\newblock \url{https://github.com/haozhg/odmd}, 2020.
\newblock Accessed: 2021-01-19.

\bibitem{streaming_tdmd}
{\sc M.~Hemati, E.~Deem, M.~Williams, C.~W. Rowley, and L.~N. Cattafesta}, {\em
  Improving separation control with noise-robust variants of dynamic mode
  decomposition}, in 54th AIAA Aerospace Sciences Meeting, 2016, p.~1103.

\bibitem{tdmd}
{\sc M.~S. Hemati, C.~W. Rowley, E.~A. Deem, and L.~N. Cattafesta}, {\em
  De-biasing the dynamic mode decomposition for applied koopman spectral
  analysis of noisy datasets}, Theoretical and Computational Fluid Dynamics, 31
  (2017), p.~349–368.

\bibitem{streaming_dmd}
{\sc M.~S. Hemati, M.~O. Williams, and C.~W. Rowley}, {\em Dynamic mode
  decomposition for large and streaming datasets}, Physics of Fluids, 26
  (2014), p.~111701.

\bibitem{kf_state_param}
{\sc R.~Isermann and M.~M{\"u}nchhof}, {\em State and Parameter Estimation by
  Kalman Filtering}, Springer Berlin Heidelberg, Berlin, Heidelberg, 2011,
  pp.~539--551.

\bibitem{tu_review}
{\sc T.~Jonathan~H., R.~Clarence~W., L.~Dirk~M., B.~Steven~L., and
  K.~J.~Nathan}, {\em On dynamic mode decomposition: Theory and applications},
  Journal of Computational Dynamics, 1 (2014), p.~391.

\bibitem{sparse_dmd}
{\sc M.~R. Jovanović, P.~J. Schmid, and J.~W. Nichols}, {\em
  Sparsity-promoting dynamic mode decomposition}, Physics of Fluids, 26 (2014),
  p.~024103.

\bibitem{kf_og}
{\sc R.~E. Kalman}, {\em {A New Approach to Linear Filtering and Prediction
  Problems}}, Journal of Basic Engineering, 82 (1960), pp.~35--45.

\bibitem{ili_mech_or_stat}
{\sc S.~Kandula, T.~Yamana, S.~Pei, W.~Yang, H.~Morita, and J.~Shaman}, {\em
  Evaluation of mechanistic and statistical methods in forecasting
  influenza-like illness}, Journal of The Royal Society Interface, 15 (2018),
  p.~20180174.

\bibitem{data_assimilation_book}
{\sc K.~J.~H. Law, A.~M. Stuart, and K.~C. Zygalakis}, {\em Data assimilation:
  A mathematical introduction}, 2015.

\bibitem{ili_mechanistic_benefits}
{\sc J.~Lessler and D.~Cummings}, {\em Mechanistic models of infectious disease
  and their impact on public health}, American Journal of Epidemiology, 183
  (2016), p.~kww021.

\bibitem{da_use_ecology}
{\sc Y.~Luo, K.~Ogle, C.~Tucker, S.~Fei, C.~Gao, S.~LaDeau, J.~S. Clark, and
  D.~S. Schimel}, {\em Ecological forecasting and data assimilation in a
  data-rich era}, Ecological Applications, 21 (2011), pp.~1429--1442.

\bibitem{dmd_autoencoder}
{\sc B.~Lusch, J.~N. Kutz, and S.~L. Brunton}, {\em Deep learning for universal
  linear embeddings of nonlinear dynamics}, Nature Communications, 9 (2018).

\bibitem{enkf_tutorial}
{\sc J.~Mandel}, {\em A brief tutorial on the ensemble kalman filter}, 2009.

\bibitem{finance_use}
{\sc J.~Mann and J.~N. Kutz}, {\em Dynamic mode decomposition for financial
  trading strategies}, 2015.

\bibitem{kernel_dmd}
{\sc I.~G.~K. "Matthew O.~Williams", "Clarence W.~Rowley"}, {\em A kernel-based
  method for data-driven koopman spectral analysis}, Journal of Computational
  Dynamics, 2 (2015), p.~247.

\bibitem{kfdmd}
{\sc T.~Nonomura, H.~Shibata, and R.~Takaki}, {\em Dynamic mode decomposition
  using a kalman filter for parameter estimation}, AIP Advances, 8 (2018),
  p.~105106.

\bibitem{ekfdmd}
{\sc T.~Nonomura, H.~Shibata, and R.~Takaki}, {\em Extended-kalman-filter-based
  dynamic mode decomposition for simultaneous system identification and
  denoising}, PloS one, 14 (2019), p.~e0209836.

\bibitem{ili_modelling_review1}
{\sc E.~O. Nsoesie, J.~S. Brownstein, N.~Ramakrishnan, and M.~V. Marathe}, {\em
  A systematic review of studies on forecasting the dynamics of influenza
  outbreaks}, Influenza and Other Respiratory Viruses, 8 (2014), pp.~309--316.

\bibitem{ili_mechanistic2}
{\sc D.~Osthus, K.~S. Hickmann, P.~C. Caragea, D.~Higdon, and S.~Y.~D. Valle},
  {\em {Forecasting seasonal influenza with a state-space SIR model}}, The
  Annals of Applied Statistics, 11 (2017), pp.~202 -- 224.

\bibitem{epidemiology_use}
{\sc J.~Proctor and P.~Welkhoff}, {\em Discovering dynamic patterns from
  infectious disease data using dynamic mode decomposition}, International
  health, 7 (2015), pp.~139--45.

\bibitem{ilinet_model_comparison}
{\sc N.~G. Reich, L.~C. Brooks, S.~J. Fox, S.~Kandula, C.~J. McGowan, E.~Moore,
  D.~Osthus, E.~L. Ray, A.~Tushar, T.~K. Yamana, M.~Biggerstaff, M.~A.
  Johansson, R.~Rosenfeld, and J.~Shaman}, {\em A collaborative multiyear,
  multimodel assessment of seasonal influenza forecasting in the united
  states}, Proceedings of the National Academy of Sciences, 116 (2019),
  pp.~3146--3154.

\bibitem{da_use_earth_sciences}
{\sc R.~H. Reichle}, {\em Data assimilation methods in the earth sciences},
  Advances in water resources, 31 (2008), pp.~1411--1418.

\bibitem{enkf_flexible}
{\sc R.~H. Reichle, J.~P. Walker, R.~D. Koster, and P.~R. Houser}, {\em
  Extended versus ensemble kalman filtering for land data assimilation},
  Journal of Hydrometeorology, 3 (01 Dec. 2002), pp.~728 -- 740.

\bibitem{ekf}
{\sc M.~Ribeiro and I.~Ribeiro}, {\em Kalman and extended kalman filters:
  Concept, derivation and properties}, 04 2004.

\bibitem{official_intro}
{\sc P.~{Schmid} and J.~{Sesterhenn}}, {\em {Dynamic Mode Decomposition of
  numerical and experimental data}}, in APS Division of Fluid Dynamics Meeting
  Abstracts, vol.~61 of APS Meeting Abstracts, Nov. 2008, p.~MR.007.

\bibitem{kernel_density_estimation}
{\sc S.~J. Sheather}, {\em Density estimation}, Statistical Science, 19 (2004),
  pp.~588--597.

\bibitem{silverman_density_estimation}
{\sc B.~W. Silverman}, {\em Density estimation for statistics and data
  analysis}, Routledge, 2018.

\bibitem{particle_filter_high_dims}
{\sc C.~Snyder, T.~Bengtsson, P.~Bickel, and J.~L. Anderson}, {\em Obstacles to
  high-dimensional particle filtering}, Monthly Weather Review, 136 (2008),
  pp.~4629--4640.

\bibitem{r0}
{\sc P.~Van~den Driessche}, {\em Reproduction numbers of infectious disease
  models}, Infectious Disease Modelling, 2 (2017), pp.~288--303.

\bibitem{ukf}
{\sc E.~A. {Wan} and R.~{Van Der Merwe}}, {\em The unscented kalman filter for
  nonlinear estimation}, in Proceedings of the IEEE 2000 Adaptive Systems for
  Signal Processing, Communications, and Control Symposium (Cat. No.00EX373),
  2000, pp.~153--158.

\bibitem{ili_statistical1}
{\sc Z.~Wang, P.~Chakraborty, S.~R. Mekaru, J.~S. Brownstein, J.~Ye, and
  N.~Ramakrishnan}, {\em Dynamic poisson autoregression for
  influenza-like-illness case count prediction}, in Proceedings of the 21th ACM
  SIGKDD International Conference on Knowledge Discovery and Data Mining, KDD
  '15, New York, NY, USA, 2015, Association for Computing Machinery,
  p.~1285–1294.

\bibitem{extended_dmd}
{\sc M.~O. Williams, I.~G. Kevrekidis, and C.~W. Rowley}, {\em A data–driven
  approximation of the koopman operator: Extending dynamic mode decomposition},
  Journal of Nonlinear Science, 25 (2015), p.~1307–1346.

\bibitem{online_dmd}
{\sc H.~Zhang, C.~W. Rowley, E.~A. Deem, and L.~N. Cattafesta}, {\em Online
  dynamic mode decomposition for time-varying systems}, 2017.

\end{thebibliography}

\end{document}